# Equations over valued groups

by Vincent Bagayoko

IMJ-PRG (Paris)

*Email:* bagayoko@imj-prg.fr

**Abstract**

We study groups, exponential groups and ordered groups equipped with valuations. We investigate algebraic and topological features of such valued structures, and apply our findings in order to solve regular equations over groups using simple valuation theoretic arguments.

# Introduction

Given a first-order theory $T$ of groups, exponential groups or ordered groups, a model $\mathcal{M}=(M,\ldots)$ of $T$, a unary term $t(y)$ in its language with parameters in $M$, the existence of an extension $\mathcal{N}\vDash T$ of $\mathcal{M}$ in which the equation

$$t(y)=1 \tag{1}$$

has a solution is a difficult problem [33]. Yet understanding the combinatorial or geometric properties of sets of solutions in $\mathcal{N}$ is preliminary in studying quantifier elimination over $T$ or existentially closed models of $T$. Suppose that $T$ is a class of pure groups. In that case $t(y)$ may be identified canonically with a word $t(y)=g_1\,y^{\alpha_1}\cdots g_n\,y^{\alpha_n}$ in the free product $\mathcal{G}*\mathbb{Z}$, where $\alpha_1,\ldots,\alpha_n\in\mathbb{Z}$ and $g_1,\ldots,g_n\in\mathcal{G}$ (imposing that $g_i\neq 1$ and $\alpha_{n-i+1}\neq 0$ for $i>1$). Set $\alpha(t):=\alpha_1+\cdots+\alpha_n\in\mathbb{Z}$. One usually distinguishes between *regular* terms, for which $\alpha(t)\neq 0$, and *singular* ones for which $\alpha(t)=0$. In the case of $A$-groups [26, 27] over a ring $A$, the $\mathbb{Z}$-valued morphism $\alpha : t\mapsto\alpha(t)$ can be extended to an $A$-valued one (see Section 7.1). We call regular those terms $t$ in the language of $A$-groups for which $\alpha(t)\neq 0$. In this paper, we use valuations on groups, exponential groups and ordered groups in order to study regular terms over such structures. This allows us to prove, for instance:

**Theorem 1.** (consequence of Theorem 7.7) *Suppose that $C$ is a field of characteristic $0$, that $A$ is a commutative divisible domain of characteristic $0$, and that $F$ is an ordered field. Let $\mathcal{G}$ be among the following exponential groups:*

a) *the $C$-group of automorphisms of $C(\!(t)\!)$ that fix the leading term of series.*

b) *the pro-nilpotent completion of a free $A$-group.*

c) *the $F$-group, under composition, of Hahn series $x+\delta, F\delta<x$ with coefficients and exponents in $F$.*

*Then for each regular term $t(y)$ over $\mathcal{G}$, there is a unique $f\in\mathcal{G}$ with $t(f)=1$.*

Let us review the difficulties of the general unary equation problem. If $T$ is a theory of Abelian groups, then $t(y)=t(1)\,y^{\alpha(t)}$ modulo $T$, so the unary equation problem is that of finding roots. In the torsion-free case, it can be conveniently solved by showing that the class of models of $T$ is preserved under taking divisible closures. In general however, the existence of solutions to even regular unary equations is more intricate: it is open whether all unary equations can be solved in extensions [33, Conjecture 2.6]. Certain specific cases do have a positive solution. In particular, divisibility equations $g\,y^{\alpha}=1$, can always be solved [28]. As for the singular case, conjugacy equations $y^{-1}\,f^{-1}\,y\,g=1$ for $f,g\in\mathcal{G}$ can always be solved using HNN-extensions [19].

The situation becomes direr if one considers totally ordered groups instead of pure groups. There, the existence of divisible extensions fails [10] (see also [13]), and until recently [5], it was an open question whether there exists even one non-trivial ordered group in which all positive elements are conjugate [11, Problem 3.31]. So far, the unary equation problem over exponential groups has received little attention. On the side of enlarging $T$ in the language of pure groups, we have the same obstructions, for orderable groups (see [17]), as in the ordered case, while various versions of the unary equation problem over torsion-free groups or solvable groups are in general open (see [33, Conjectures 2.7 and 2.17]).

We introduced [8] an elementary class of ordered groups called *growth order groups* (see Section 4.4) over which we expect that the unary equation problem should be more traceable. Prominent examples are groups under composition of germs of functions definable in o-minimal structures [31], and groups under formal composition laws of generalised power series, such as transseries [16, 15]. Growth order groups come equipped with a canonical definable valuation [35], that contains a lot of information on solutions of unary equations. Those that come from fields of generalised power series often have a natural non-trivial structure of exponential group, which is compatible with the valuation. In general, solving unary equations over growth order groups requires considering intermediate extensions that fail to be growth order groups in general, though they naturally inherit a well-behaved valuation. Furthermore, some of the benefits of working with growth order groups come in large part from the properties of their canonical valuation, irrespective of the ordering. It thus appears useful to study in their own right groups with abstract valuations that behave like the canonical valuation on a growth order group, and that encompass extensions of growth order groups that are involved in solving unary equations.

Accordingly, we study a class of valued groups [35], called *c-valued groups*, where the valuation $v$ is compatible with both commutation (i.e. it is invariant on centralisers) and commutators (i.e. commutators of elements with the same valuation decrease in valuation). Our goals are to show that (ordered) c-valued groups form a natural class of expansions of (ordered) groups over which it is easier to study unary equations (and inequalities), and to give examples of growth order groups in which regular equations can be solved.

In commutative algebra, dating to the Newton-Puiseux polygon method (see [4, Section 3.8]), valuations have proven very useful for studying algebraic and differential-algebraic equations. Let $K$ be a field or an H-field [1, 2, 4], and let $t(y)$ be a unary term over $K$ in the corresponding first-order language. If $K$ is endowed with a valuation $v: K^{\times} \longrightarrow \Gamma$, then the possible valuations of a solution of $t(y)=0$ in $K$ can be understood in terms of valuation theoretic characteristics of $t$. Fixing the valuation of a solution, the corresponding possible residues must also satisfy algebraic equations in the residue field $Kv$, which in applications is usually better understood. Furthermore, the existence of a solution can be inferred from closure properties of $\Gamma$ (i.e. divisibility of the value group and H-closedness of the asymptotic couple [4, Section 9.9]) and $Kv$ (i.e. algebraic closedness or real closedness) together with algebraic completeness properties (i.e. Henselianity [32, Section 3.2] and Newtonianity [4, Chapter 14]) of $(K, v)$. Valuations on groups can be used in the same fashion. Given a c-valuation $v: \mathcal{G} \longrightarrow \Gamma$ on a(n ordered) group $\mathcal{G}$, to each $\rho \in \Gamma$ corresponds a(n ordered) Abelian group $\mathcal{C}_{\rho}$ of *residues*. For regular terms $t(y)$, the existence of a solution to $t(y)=1$ can be inferred from the existence of solutions to equations in $\mathcal{C}_{\rho}$, i.e. that it be divisible, along with a completeness condition called spherical completeness. Moreover, the valuation and residue of a solution, as well as the behaviour of the function $\mathcal{G} \ni f \mapsto t(f)$ can be inferred from simple valuation-theoretic characteristics of $t(y)$. This also applies to valued exponential groups. This is, roughly speaking, the content of our main result Theorem 2, which we describe hereafter.

We give conventions on groups and orderings in Section 1. Section 2 introduces our notions of valuations and c-valuations on groups. We give examples, derive group theoretic consequences of the existence of a c-valuation on a group, consider c-valued quotients, and define the Abelian residue groups $\mathcal{C}_\rho$, $\rho \in v(\mathcal{G})$. Sections 3 and 4 establish for exponential groups and ordered groups respectively what we established for groups in Section 2.

Section 5 recalls standard notions of valuative topology and uniform structure on valued groups (see [35]), with emphasis on the case of c-valued groups and $A$-valued $A$-groups. In particular, we show (Theorem 5.11) that the completion of a valued group with respect to its uniform structure has a natural structure of valued group. Similar results hold for $A$-valued $A$-groups (Proposition 5.12), as well as for $c$-valued groups and ordered c-valued groups under a further condition (see Definition 2.27). In Section 5.3, we study the notion of pseudo-Cauchy sequence and pseudo-limit as in [35]. Spherically complete valued groups are valued groups in which pseudo-Cauchy sequences have pseudo-limits.

In Section 6, we study a class of valued groups called *nearly Abelian* valued groups, in which commutators are small relative to their arguments. Lastly, we study in Section 7 regular equations over nearly Abelian valued exponential groups. Our main results are derived from a technical description (Lemma 7.5) of the possible valuations of evaluations of regular unary terms. They can summed up as follows:

**Theorem 2.** (consequence of Corollary 7.6 and Theorems 7.9 and 7.7) *Let $A$ be a ring of characteristic* 0. *Let $(\mathcal{G}, \cdot, 1, v)$ be an $A$-torsion-free nearly Abelian $A$-valued $A$-group and let $t(y)$ be a regular unary term over $\mathcal{G}$.*

*a) There is at most one solution of $t(y) = 1$ in $\mathcal{G}$.*

*b) If $A = \mathbb{Z}$ and $\mathcal{G}$ is an ordered valued group, then the function $\mathcal{G} \longrightarrow \mathcal{G}; f \mapsto t(f)$ is strictly monotonous.*

*c) If $\mathcal{G}$ is spherically complete and each residue $A$-module $\mathcal{C}_\rho$, $\rho \in v(\mathcal{G})$ is divisible, then there is a unique solution of $t(y) = 1$ in $\mathcal{G}$.*

The problem of solvability of regular unary equations over torsion-free nearly Abelian valued groups is thus reduced to showing how to embed them into spherically complete and divisible ones, which we will do in future work. We also deduce an extension to exponential groups (Theorem 7.10) of an old result of Smel'kin [37] on unimodular equations over torsion-free nilpotent groups.

# 1 Preliminaries

## 1.1 Group theoretic conventions

Let $(\mathcal{G}, \cdot, 1)$ be a group. The *commutator* of an $(f, g) \in \mathcal{G} \times \mathcal{G}$ is the element

$$[f, g] := f^{-1} g^{-1} f g \in \mathcal{G}.$$

The *centraliser* of an $f \in \mathcal{G}$ is the subgroup

$$\mathcal{C}(f) := \{g \in \mathcal{G} : [f, g] = 1\} = \{g \in \mathcal{G} : f g = g f\}.$$

of $\mathcal{G}$. We write

$$\mathcal{G}^{\neq} := \{f \in \mathcal{G} : f \neq 1\}.$$

Our first-order language of groups $\mathcal{L}_{\mathrm{g}} = \langle \cdot, 1, \mathrm{Inv} \rangle$ consists of a binary function symbol $\cdot$ interpreted as the group law, a constant symbol 1 interpreted as the identity, and a unary function symbol Inv interpreted as the inverse map. Note that $\mathcal{L}_{\mathrm{g}}$-substructures of groups are simply subgroups. A group $(\mathcal{G}, \cdot, 1)$ is said *divisible* if for all $g \in \mathcal{G}$ and all $n \in \mathbb{Z} \setminus \{0\}$, there is an $f \in \mathcal{G}$ with $f^n = g$.

Given a group $(\mathcal{G}, \cdot, 1)$, the *inverse group* $\check{\mathcal{G}}$ is the group $(\mathcal{G}, \check{\cdot}, 1)$ where $f \,\check{\cdot}\, g := g \cdot f$ for all $f, g \in \mathcal{G}$. The inverse map on $\mathcal{G}$ is an isomorphism between $\mathcal{G}$ and $\check{\mathcal{G}}$.

### 1.2 Orderings and quasi-orderings

An *ordering* is an anti-reflexive and transitive binary relation. A *quasi-ordering* is a reflexive and transitive binary relation. Given a quasi-ordering $\preccurlyeq$ on a set $X$, we usually write $\asymp$ for the corresponding equivalence relation where $f \asymp g \Longleftrightarrow f \preccurlyeq g \wedge g \preccurlyeq f$ for all $f$, $g \in X$. We write $\prec$ for the corresponding ordering where $f \prec g \Longleftrightarrow (f \preccurlyeq g \wedge g \not\preccurlyeq f)$ for all $f, g \in X$. If $\preccurlyeq$ is total, then we have $f \prec g \Longleftrightarrow g \not\preccurlyeq f$ for all $f, g \in X$.

Given two totally ordered sets $(\Gamma_1, <_1)$ and $(\Gamma_2, <)$, the disjoint sum $\Gamma_1 \amalg \Gamma_2$ is the union $(\Gamma_1 \times \{0\}) \sqcup (\Gamma_2 \times \{1\})$ with the total ordering $<$ extending $<_1$ on $\Gamma_1 \times \{0\}$ and $<_2$ on $\Gamma_2 \times \{1\}$ and with $\Gamma_1 \times \{0\} < \Gamma_2 \times \{1\}$.

## 2 Valuations on groups

### 2.1 Groups with dominance relations

**Definition 2.1.** *A* **dominance relation** *on a group* $(\mathcal{G}, \cdot, 1)$ *is total quasi-ordering* $\preccurlyeq$ *on* $\mathcal{G}$ *such for all* $f, g, h \in \mathcal{G}$*, we have*

**D1.** $f \neq 1 \longrightarrow 1 \prec f$.

**D2.** $f g \preccurlyeq f$ *or* $f g \preccurlyeq g$.

**D3.** $f \preccurlyeq g \longrightarrow h f h^{-1} \preccurlyeq h g h^{-1}$.

**D4.** $f \asymp f^{-1}$.

*A* **c-dominance relation** *on a group* $(\mathcal{G}, \cdot, 1)$ *is a dominance relation* $\preccurlyeq$ *on* $\mathcal{G}$ *such for all* $f, g, h \in \mathcal{G}$*, we have*

**D5.** $(g \in \mathcal{C}(f) \wedge f \neq 1 \wedge g \neq 1) \longrightarrow g \asymp f$.

**D6.** $(f \asymp g \wedge f \neq 1 \wedge g \neq 1) \longrightarrow [f, g] \prec f$.

*A* **(c-)valued group** *is a group* $(\mathcal{G}, \cdot, 1)$ *together with a (c-)dominance relation.*

We see a valued group $(\mathcal{G}, \cdot, 1, \preccurlyeq)$ as an $\mathcal{L}_{\mathrm{vg}}$-structure where $\mathcal{L}_{\mathrm{vg}}$ is the expansion of $\mathcal{L}_{\mathrm{g}}$ with a binary relation symbol $\preccurlyeq$. The second part of the axiom **D3** means that conjugating with an element in $\mathcal{G}$ is an automorphism of $\mathcal{L}_{\mathrm{vg}}$-structure. We say that an $\mathcal{L}_{\mathrm{vg}}$-sentence $\psi$ is *preserved under coarsenings* if for groups $(\mathcal{G}, \cdot, 1)$ and all total quasi-orderings $\preccurlyeq_0, \preccurlyeq_1$ such that $\preccurlyeq_1$ is coarser than $\preccurlyeq_0$, we have $(\mathcal{G}, \cdot, 1, \preccurlyeq_1) \vDash \psi$ if $(\mathcal{G}, \cdot, 1, \preccurlyeq_0) \vDash \psi$.

**Remark 2.2.** The axioms **D1**, **D2**, **D4** and **D5** are preserved under coarsenings.

If $(\mathcal{G}, \cdot, 1, \preccurlyeq)$ is (c-)valued, then we have a quotient set $\mathrm{val}(\mathcal{G}) := \mathcal{G}^{\neq} / \asymp$. We usually denote the quotient map by val. Since $\preccurlyeq$ is a total quasi-ordering, the set $\mathrm{val}(\mathcal{G})$ is totally ordered by $\mathrm{val}(f) < \mathrm{val}(g) \Longleftrightarrow f \prec g$. We have the following equivalent formulations for val of the axioms **D1**–**D6** for all $f, g, h \in \mathcal{G}$:

**V1.** $f \neq 1 \longrightarrow \mathrm{val}(1) < \mathrm{val}(f)$.

**V2.** $\mathrm{val}(f g) \leqslant \max(\mathrm{val}(f), \mathrm{val}(g))$.

**V3.** $\mathrm{val}(f) \leqslant \mathrm{val}(g) \longleftrightarrow \mathrm{val}(h f h^{-1}) \leqslant \mathrm{val}(h g h^{-1})$.

**V4.** $\mathrm{val}(f) = \mathrm{val}(f^{-1})$

**V5.** $(g \in \mathcal{C}(f) \wedge f \neq 1 \wedge g \neq 1) \longrightarrow \mathrm{val}(g) = \mathrm{val}(f)$.

**V6.** $(\mathrm{val}(f) = \mathrm{val}(g) \wedge f \neq 1 \wedge g \neq 1) \longrightarrow \mathrm{val}([f, g]) < \mathrm{val}(f)$.

Conversely, if $(\Gamma, <)$ is a totally ordered set and $v : \mathcal{G} \longrightarrow \Gamma$ is a surjective map satisfying **V1–V4** (resp. **V1–V6**) then we say that $v$ is a valuation (resp. *c-valuation*) on $\mathcal{G}$, and we call $\Gamma \setminus \{\min(\Gamma)\}$ the *value set* of $v$. If $v$ is a (c-)valuation, then the induced relation $f \preccurlyeq g \Longleftrightarrow v(f) \leqslant v(g)$ is a (c-)dominance relation on $\mathcal{G}$. Thus the settings of (c-)valuations and (c-)dominance relations are equivalent.

**Remark 2.3.** In [35, Definition 2.1], a valuation on a group is defined to satisfy just **V1** and **V2**. Our c-valuations are much more restrictive.

**Proposition 2.4.** *Let $(\mathcal{G}, \cdot, 1, \preccurlyeq)$ be a valued group and $f, g \in \mathcal{G}$ with $f \prec g$. Then $f g \asymp g f \asymp g$.*

**Proof.** We have $f g, g f \preccurlyeq g$ by **D2**. Suppose for contradiction that $f g \prec g$. We have $g = f^{-1}(f g) \preccurlyeq f g$ or $g \preccurlyeq f^{-1}$ by **D2**, so we muse have $g \preccurlyeq f^{-1}$. But $f^{-1} \asymp f$ by **D4**, so $g \preccurlyeq f$: a contradiction. Symmetric arguments show that $g f \not\prec g$, hence the result. $\square$

The following immediate consequence of **D3** is sometimes useful.

**Lemma 2.5.** *Let $\mathcal{G}$ be a valued group. For $f, g \in \mathcal{G}$ with $f \prec g$, we have*

$$g f g^{-1} \prec g.$$

**Lemma 2.6.** *Let $\mathcal{G}$ be a valued group. Let $f, g \in \mathcal{G}$ with $f \not\asymp g$. We have* $\mathrm{val}([f, g]) < \max(\mathrm{val}(f), \mathrm{val}(g))$.

**Proof.** Since $[f, g] = [g, f]^{-1} \asymp [g, f]$ by **D4**, we may assume that $g \prec f$, i.e. $g^{-1} \prec f^{-1}$ by **D4**. By Lemma 2.5, we have $f^{-1} g^{-1} f \prec f$, whence $f^{-1} g^{-1} f g \prec f$ by Proposition 2.4. $\square$

Let $(\mathcal{G}, \cdot, 1, \preccurlyeq)$ be a valued group and let $\rho \in \mathrm{val}(\mathcal{G})$. By **D3**, for all $f, g \in \mathcal{G}$ with $\rho = \mathrm{val}(g)$ the value $v(f g f^{-1})$ only depends on $f$ and $\rho$. We write $f \rho f^{-1} := \mathrm{val}(f g f^{-1})$.

A homomorphism (resp. an embedding) of valued groups is a homomorphism (resp. an embedding) of $\mathcal{L}_{\mathrm{vg}}$-structures. If $\mathcal{G}$ and $\mathcal{H}$ are two valued groups with respective valuations and $\Phi : \mathcal{H} \longrightarrow \mathcal{G}$ is a homomorphism (resp. an embedding) of valued groups, then it induces a nondecreasing (resp. strictly increasing map)

$$\begin{aligned} \mathrm{val}_\Phi : \mathrm{val}(\mathcal{G}) &\longrightarrow \mathrm{val}(\mathcal{H}) \\ \mathrm{val}(f) &\longmapsto \mathrm{val}(\Phi(f)). \end{aligned}$$

## 2.2 Examples of c-valued groups

We give a few examples of c-valued groups.

**Example 2.7.** The theory of c-valued Abelian groups is trivial in the sense that every Abelian group has a unique c-dominance relation $\preccurlyeq$ given by $f \preccurlyeq g$ for all $f \in \mathcal{G}$ and $g \in \mathcal{G}^{\neq}$.

**Proposition 2.8.** *A subgroup of a c-valued group is c-valued for the induced quasi-ordering.*

**Proof.** The axioms of c-valued groups can be formulated as universal sentences in $\mathcal{L}_{\mathrm{vg}}$. Since $\mathcal{L}_{\mathrm{vg}}$-substructures are just subgroups, the result follows. □

**Remark 2.9.** Let $(\mathcal{G}, v), (\mathcal{H}, w)$ be two non-trivial valued groups. Define a map $u: \mathcal{H}\times\mathcal{G}\longrightarrow w(\mathcal{H})\amalg v(\mathcal{G}^{\neq})$ by $u(h,1)=w(h)$ and $u(h,g)=v(g)$ for all $h\in\mathcal{H}$ and $g\in\mathcal{G}^{\neq}$. Then $u$ is a valuation on $\mathcal{H}\times\mathcal{G}$. However, it is not a c-valuation. Indeed, given $(h,g)\in\mathcal{H}^{\neq}\times\mathcal{G}^{\neq}$, the non-trivial elements $(h,1)$ and $(1,g)$ commute but have distinct valuations. So **V5** does not hold. In view of Proposition 2.8, no non-trivial product or Hahn product construction of valued groups (see [35, Section 2]) yields a c-valued group.

**Example 2.10.** Let $F$ be a field and let $V$ be a non-trivial vector space over $k$. Then the group under composition of bijective affine maps $\lambda x+v: u\mapsto \lambda u+v$ for $(\lambda,v)\in F^{\times}\times V$ is denoted $\mathrm{Aff}(V)$. This can be seen as the semidirect product $F^{\times}\rtimes V$ for the action of $F^{\times}$ on $V$ by scalar multiplication. This group is solvable of derived length 2. Indeed $[\mathrm{Aff}(V),\mathrm{Aff}(V)]\subseteq x+V$ is Abelian.

We define a valuation $v:\mathrm{Aff}(V)\longrightarrow(\{0,1,2\},>)$ by $v(\lambda x+v):=0$ if $\lambda\neq 1$, $v(x+v)=1$ if $v\neq 0$ and $v(x)=2$. It satisfies **V5** as $\mathcal{C}(x+v)=x+V$ for all $v\in V\setminus\{0\}$ whereas $\mathcal{C}(\lambda x+v)\subseteq(F^{\times}\setminus\{1\})x+V$ if $\lambda\neq 1$ (one computes centralisers in $\mathrm{Aff}(V)$ by solving simple affine equations). Since $[\mathrm{Aff}(V),\mathrm{Aff}(V)]\subseteq x+V$ and $[x+V,x+V]=\{1\}$, the axiom **V6** holds as well. Thus $(\mathrm{Aff}(V),\circ,x,v)$ is a c-valued group. Note that $\mathrm{Aff}(V)$ has torsion.

**Example 2.11.** Let $C$ be a field of characteristic zero, and let $C((t))$ be the field of formal Laurent series with coefficients in $C$. It is equipped with its canonical valuation $v: C((t))^{\times}\longrightarrow\mathbb{Z}$ (in the sense of commutative algebra) sending a non-zero series to the smallest exponent in $\mathbb{Z}$ of a non-zero coefficient. There is a well-defined composition law $\circ: C((t))\times C((t))^{\geqslant 1}\longrightarrow C((t))$ where $C((t))^{\geqslant 1}$ is the set of series with valuation $\geqslant 1$. For $f=\sum_{i\geqslant m} f_i t^i\in C((t))$ and $s\in C((t))^{\geqslant 1}$, the series $f\circ s$ is given by the valuation theoretic limit of the sequence $\left(\sum_{i\in[m,m+n]} f_i t^i\right)_{n\in\mathbb{N}}$. It is easy to see that $v(f\circ s)=v(f)\,v(s)$. Now let $C((t))^{\equiv 1}$ be the set of series with valuation 1 and leading coefficient 1. Note that $C((t))^{\equiv 1}$ is closed under composition, whence they form monoids with identity $t$. Furthermore, for all $f\in C((t))^{\equiv 1}$ and $s=t+\delta\in C((t))^{\equiv 1}$, we have

$$f\circ s=\sum_{i\in\mathbb{N}}\frac{f^{(i)}}{i!}\,\delta^i, \tag{2.1}$$

where $f^{(i)}$ denotes the $i$-th derivative of $f$ with respect to $t$. Indeed this holds for $f=t$ and the linear operator $f\mapsto\sum_{i\in\mathbb{N}}\frac{f^{(i)}}{i!}\,\delta^i$ commutes with infinite sums and finite products. Thus the function $\circ_s: C((t))^{\equiv 1}\longrightarrow C((t))^{\equiv 1};\ g\mapsto g\circ s$ is of the form $\circ_s=\mathrm{Id}+\phi$ where

$$v(\phi(f))=v(f\circ s-f)\geqslant v(f)+v(\delta)=v(f)-1+v(\delta)\geqslant v(f)-1+2>v(f)$$

for all $f\in C((t))^{\geqslant 1}$. By [3, Corollary 1.4] it has an inverse $\Psi$, which means that the series $\Psi(t)$ is the inverse of $s$ in $C((t))^{\equiv 1}$. Therefore $C((t))^{\equiv 1}$ is a group. Given $f\in C((t))^{\equiv 1}$, we write $f^{\mathrm{inv}}$ for its inverse for the composition law.

We define a c-valuation val on $C((t))^{\equiv 1}$ with value set $(\omega\setminus\{0,1\},>)$ as follows

$$\forall f\in C((t))^{\equiv 1}\setminus\{t\}, \mathrm{val}(f):=v(f-t).$$

We let the reader check that **V1** and **V2** hold for val. For all $f,g\in C((t))^{\equiv 1}$, we have

$$f\circ g=t+(g-t)+(f-t)+(f-t)'\,(g-t)+\cdots \tag{2.2}$$

where $v((f-t)'(g-t)+\cdots)$ is strictly larger than both $v(f-t)$ and $v(g-t)$. From these observations, we deduce that:

i. $\mathrm{val}(f\circ g)\leqslant\max(v(g-t),v(f-t))=\max(\mathrm{val}(f),\mathrm{val}(g))$.

ii. $f^{\mathrm{inv}}=t-(f-t)+\varepsilon$ for an $\varepsilon$ with $v(\varepsilon)>v(f-t)$. Hence:

iii. $f\circ g\circ f^{\mathrm{inv}}=t+(g-t)+\delta$ for a $\delta$ with $v(\delta)>v(g-t),v(f-t)$, whence:

iv. $\mathrm{val}(f\circ g\circ f^{\mathrm{inv}})=\mathrm{val}(g)$.

v. If $f$ and $g$ are non-trivial and commute, then $(f-t)'(g-t)+\cdots$ and $(g-t)'(f-t)+\cdots$ are equal, so they have the same leading coefficient. Writing $a,b\in C^{\times}$ for the leading coefficient of $f-t$ and $g-t$ respectively, this gives the identity $\mathrm{val}(f)\,a\,b=\mathrm{val}(g)\,a\,b$, whence $\mathrm{val}(f)=\mathrm{val}(g)$.

vi. In view of iii, we have $\mathrm{val}(f\circ g\circ f^{\mathrm{inv}}\circ g^{\mathrm{inv}})<_{\Gamma}\mathrm{val}(f),\mathrm{val}(g)$ when $f,g\neq t$.

So **V2**–**V6** hold by i, iv, ii, v and vi respectively, and $(C((t))^{\equiv 1},\circ,t,\mathrm{val})$ is a c-valued group.

## 2.3 Basic properties of c-valued groups

We now give a few important properties of c-valued groups, including purely group theoretic ones. We fix a $c$-valued group $(\mathcal{G},\cdot,1,\preccurlyeq)$.

**Proposition 2.12.** *For each $f\in\mathcal{G}^{\neq}$, the subgroup $\mathcal{C}(f)$ is Abelian.*

**Proof.** Let $g,h\in\mathcal{C}(f)^{\neq}$. So $g\asymp h$ by **D5**. Then **D6** gives $[g,h]\prec f$, but $[g,h]\in\mathcal{C}(f)$ so **D5** entails that $[g,h]=1$. □

So c-valued groups are CT-groups as per [38].

**Corollary 2.13.** *No non-Abelian nilpotent group can be c-valued.*

**Lemma 2.14.** *For all $f,g\in\mathcal{G}$ and $n\in\mathbb{Z}$ with $f^n\neq 1$, we have $[f^n,g]=1\Longrightarrow[f,g]=1$.*

**Proof.** We may assume that $f,g\neq 1$ and that $n>0$. Write $\varepsilon=[f,g]$ and suppose for contradiction that $\varepsilon\neq 1$. Note that $f^n\asymp f\asymp g$ by **D5**, so $\varepsilon\prec g$ by **D6**. We have $f^{-1}g^{-1}f=\varepsilon\,g^{-1}$, so

$$1=f^{-n}\,g\,f^n=f^{-(n-1)}\,\varepsilon\,f^{n-1}\,f^{-(n-1)}\,g^{-1}\,f^{n-1}.$$

But $\varepsilon\prec g\asymp g^{-1}$ so $f^{-(n-1)}\,\varepsilon\,f^{n-1}\prec f^{-(n-1)}\,g^{-1}\,f^{n-1}$ by **D3**. Thus $1\asymp f^{-(n-1)}\,g^{-1}\,f^{n-1}$ by Proposition 2.4. But we have $f^{-(n-1)}\,g^{-1}\,f^{n-1}\asymp g^{-1}\not\asymp 1$ by **D1**: a contradiction. □

**Proposition 2.15.** *Suppose that $(\mathcal{G},\cdot,1,\preccurlyeq)$ has value set $(\lambda,>)$ for some ordinal $\lambda$. Then $\mathcal{G}$ is $\mathbb{Z}$-hypoabelian with $\mathcal{G}^{(\lambda)}=\{1\}$.*

**Proof.** We may assume that $\lambda>0$, for otherwise $\mathcal{G}$ is trivial. Identifying $(\lambda,>)$ and $\mathrm{val}(\mathcal{G})$, we prove by induction on an ordinal $\alpha$ that for all $\gamma\leqslant\alpha$ and $f\in\mathcal{G}_{\gamma}$, we have $\mathrm{val}(f)\geqslant\gamma$. It then follows that $\mathcal{G}_{\lambda}=\{1\}$, since any $f\in\mathcal{G}^{\neq}$ has valuation $\mathrm{val}(f)\in\lambda$. Suppose that the result holds for all $\beta<\alpha$. If $\alpha$ is a limit, then the result follows immediately by induction. Suppose that $\alpha=\beta+1$ is a successor. Given $f\in\mathcal{G}_{\alpha}$, there are $g,h\in\mathcal{G}_{\beta}$ with $f=[g,h]$. The induction hypothesis gives $\max\limits_{(\lambda,>)}(\mathrm{val}(g),\mathrm{val}(h))\geqslant\beta$, so $\mathrm{val}(f)>\beta$ by Lemma 2.6. But $\mathrm{val}(f)\in\lambda$ so $\mathrm{val}(f)\geqslant\beta+1=\alpha$. The result follows by induction. □

In particular, c-valued groups of finite value set (e.g. finite c-valued groups) are solvable as groups. In view of Proposition 2.8, we obtain:

**Corollary 2.16.** *Every finite subgroup of a c-valued group is solvable.*

## 2.4 Quotients and group extensions

Let $(\mathcal{G},\cdot,1,\preccurlyeq)$ be a valued groupp. A subset $\mathcal{H}\subseteq\mathcal{G}$ is said $\preccurlyeq$-initial if it is downward closed for $\preccurlyeq$, i.e. if we have

$$\forall f,g\in\mathcal{G}\,((f\preccurlyeq g\wedge g\in\mathcal{H})\Longrightarrow f\in\mathcal{H}).$$

**Lemma 2.17.** *Any non-empty $\preccurlyeq$-initial subset of $\mathcal{G}$ is a subgroup.*

**Proof.** This follows from **D1**, **D2** and **D4**. □

**Proposition 2.18.** *Let $\mathcal{N}\trianglelefteq\mathcal{G}$ be a normal and $\preccurlyeq$-initial subgroup. Then $\mathcal{G}/\mathcal{N}$ is valued for the dominance relation $f\mathcal{N}\preccurlyeq g\mathcal{N}\Longleftrightarrow(f\in\mathcal{N}\vee f\preccurlyeq g)$. If $\mathcal{G}$ satisfies* **D6***, then so does $\mathcal{G}/\mathcal{N}$. If moreover the following holds*

$$\forall f,g\in\mathcal{G}\setminus\mathcal{N},[f,g]\in\mathcal{N}\Longrightarrow f\asymp g, \tag{2.3}$$

*then $\mathcal{G}/\mathcal{N}$ is satisfies* **D5***.*

**Proof.** The relation $\preccurlyeq$ is clearly reflexive and total. Let $f,g,h\in\mathcal{G}$ with $f\mathcal{N}\preccurlyeq g\mathcal{N}$ and $g\mathcal{N}\preccurlyeq h\mathcal{N}$. If $f\in\mathcal{N}$ then $f\mathcal{N}\preccurlyeq h\mathcal{N}$. Otherwise $f\preccurlyeq g$ so $g\notin\mathcal{N}$ by $\preccurlyeq$-initiality of $\mathcal{N}$. So $g\preccurlyeq h$, whence $f\preccurlyeq h$, so $f\mathcal{N}\preccurlyeq h\mathcal{N}$. Thus $\preccurlyeq$ is a total quasi-ordering on $\mathcal{G}/\mathcal{N}$.

The axiom **D1** holds by definition, and **D4** holds by **D4** in $\mathcal{G}$. For $f,g\in\mathcal{G}$, we either have $fg\preccurlyeq f$, in which case $fg\mathcal{N}\preccurlyeq f\mathcal{N}$, or $fg\preccurlyeq g$, in which case $fg\mathcal{N}\preccurlyeq g\mathcal{N}$. So **D2** holds. That **D3** holds follows from **D3** for $\mathcal{G}$ and the fact that $\mathcal{N}$ is normal. Suppose that **D6** holds in $\mathcal{G}$. given $f,g\in\mathcal{G}\setminus\mathcal{N}$ with $f\mathcal{N}\asymp g\mathcal{N}$, we must have $f\asymp g$, so $[f,g]\prec f$. So $[f,g]\mathcal{N}=[f\mathcal{N},g\mathcal{N}]\preccurlyeq f\mathcal{N}$. We cannot have $f\mathcal{N}\preccurlyeq[f,g]\mathcal{N}$ since $f\notin\mathcal{N}$ and $[f,g]\prec f$. So $[f\mathcal{N},g\mathcal{N}]\prec f\mathcal{N}$, i.e. **D6** holds in $\mathcal{G}/\mathcal{N}$.

Lastly, suppose that (2.3) holds. Let $f,g\in\mathcal{G}\setminus\mathcal{N}$ with $f\mathcal{N}\in\mathcal{C}(g\mathcal{N})$. So $f,g\in\mathcal{G}\setminus\mathcal{N}$ and that $[f,g]\in\mathcal{N}$. Thus $f\asymp g$ in $\mathcal{G}$, whence $f\mathcal{N}\asymp g\mathcal{N}$ in $\mathcal{N}$. Therefore **D5** holds. □

Let $\rho\in\operatorname{val}(\mathcal{G})$ and write

$$\mathcal{G}_{\leqslant\rho}:=\{f\in\mathcal{G}:\operatorname{val}(f)\leqslant\rho\}\qquad\text{and}\qquad\mathcal{G}_{<\rho}:=\{f\in\mathcal{G}:\operatorname{val}(f)<\rho\}.$$

Then $\mathcal{G}_{\leqslant\rho}$ is a valued group for the induced quasi-ordering, and $\mathcal{G}_{<\rho}$ is a non-empty $\preccurlyeq$-initial subset.

**Proposition 2.19.** *The quotient $\mathcal{G}_{\leqslant\rho}/\mathcal{G}_\rho$ is a valued group satisfying* **D5** *for the quotient quasi-ordering, and it is c-valued if $\mathcal{G}$ is c-valued. Moreover, we have $\operatorname{val}(\mathcal{G}_{\leqslant\rho}/\mathcal{G}_\rho)=1$.*

**Proof.** Let $f\in\mathcal{G}_{\leqslant\rho}\setminus\mathcal{G}_{<\rho}$ and $g\in\mathcal{G}_{<\rho}$. So $\rho=\operatorname{val}(f)$. We have $fgf^{-1}\in\mathcal{G}_{<\rho}$ by Lemma 2.5, so $\mathcal{G}_{<\rho}$ is a normal subgroup of $\mathcal{G}_{\leqslant\rho}$. By definition, the condition (2.3) in Proposition 2.18 is trivially satisfied. Therefore $\mathcal{G}_{\leqslant\rho}/\mathcal{G}_\rho$ is a valued group satisfying **D5** for the quotient quasi-ordering. It is clear that $\operatorname{val}(\mathcal{G}_{\leqslant\rho}/\mathcal{G}_\rho)=1$. □

**Proposition 2.20.** *Let $(\mathcal{N},\cdot,1,\preccurlyeq_{\mathcal{N}})$ and $(\mathcal{Q},\cdot,1,\preccurlyeq_{\mathcal{Q}})$ be c-valued groups, let $(\mathcal{G},\cdot,1)$ be a group and suppose we have a short exact sequence of groups*

$$0\longrightarrow\mathcal{N}\overset{\iota}{\longrightarrow}\mathcal{G}\overset{\pi}{\longrightarrow}\mathcal{Q}\longrightarrow 0,$$

*where each conjugacy map* $\mathcal{N} \longrightarrow \mathcal{N}\,;\, h \mapsto g\,h\,g^{-1}$ *for* $g \in \mathcal{G}$ *is in* $\mathrm{Aut}(\mathcal{N}, \cdot, 1, \preccurlyeq)$ *and where* $\pi(\iota(g\,h)) \asymp_{\mathcal{Q}} \pi(\iota(g))$ *for all* $g \in \mathcal{G}$ *and* $h \in \mathcal{N}$*.*

*There is a unique dominance relation* $\preccurlyeq$ *on* $\mathcal{G}$ *for which* $\iota$ *is an embedding of valued groups and* $\pi$ *is a homomorphism of valued groups, and we have* $\mathrm{val}(\mathcal{G}) = \mathrm{val}(\mathcal{N}) \amalg \mathrm{val}(\mathcal{Q})$*. Suppose that the following condition is satisfied:*

$$\forall g \in \mathcal{G} \setminus \mathcal{N}, \forall h \in \mathcal{N}^{\neq}, \forall f \in \mathcal{N}, g\,h\,g^{-1} \neq f\,h\,f^{-1}. \tag{2.4}$$

*Then* $(\mathcal{G}, \cdot, 1, \preccurlyeq)$ *is c-valued.*

**Proof.** We may assume that $\mathcal{N} \subseteq \mathcal{G}$ and that $\iota$ is the inclusion map. Write $\Gamma = (\{\mathrm{val}_{\mathcal{H}}(1)\} \cup \mathrm{val}_{\mathcal{H}}(\mathcal{H})) \amalg \mathrm{val}_{\mathcal{G}}(\mathcal{G})$. For $f, g \in \mathcal{G}$, we set $f \preccurlyeq g$ if $f, g \in \mathcal{N}$ and $f \preccurlyeq_{\mathcal{N}} g$, or $f \in \mathcal{N}$ and $g \notin \mathcal{N}$, or $f, g \notin \mathcal{N}$ and $\pi(f) \preccurlyeq_{\mathcal{Q}} \pi(g)$. This is well-defined by our assumption on $\pi \circ \iota$. Note that $\preccurlyeq$ is the unique total quasi-ordering on $\mathcal{G}$ satisfying the conditions. It is easy to check that **D1–D4** hold. For **D3**, one uses the assumption on the conjugacy action $\mathcal{G} \longrightarrow \mathrm{Aut}(\mathcal{N})$. We clearly have $\mathcal{N} \prec \mathcal{G} \setminus \mathcal{N}$ and $\mathrm{val}(\mathcal{G} \setminus \mathcal{N}) = \mathrm{val}(\mathcal{Q})$, so $\mathrm{val}(\mathcal{G}) = \mathrm{val}(\mathcal{N}) \amalg \mathrm{val}(\mathcal{Q})$. Let $g_0$, $g_1 \in \mathcal{G}^{\neq}$ with $[g_0, g_1] = 1$. So, $[\pi(g_0), \pi(g_1)] = 1$. Assume for contradiction that $\pi(g_0) = 1$ and $\pi(g_1) \neq 1$. So $g_1 \in \mathcal{G} \setminus \mathcal{N}$ and $g_0 \in \mathcal{N}$, whence $g_1\, g_0\, g_1^{-1} \neq g_0\, g_0\, g_0^{-1} = g_0$ by (2.4): a contradiction. Likewise, it cannot be that $\pi(g_1) = 1$ and $\pi(g_0) \neq 1$. So $\pi(g_0), \pi(g_1) \in \mathcal{Q}^{\neq}$, whence $g_0 \asymp_{\mathcal{Q}} g_1$ by **D5** in $\mathcal{Q}$. This means that $g_0 \asymp g_1$ in $\mathcal{G}$, so **D5** holds in $\mathcal{G}$.

Lastly, suppose that $g_0, g_1 \in \mathcal{G}^{\neq}$ and $g_0 \asymp g_1$. Suppose that $g_0 \notin \mathcal{N}$. So also $g_1 \notin \mathcal{N}$. If $[g_0, g_1] \in \mathcal{N}$, then $[g_0, g_1] \prec g_0$. If $[g_0, g_1] \notin \mathcal{N}$, then since $\pi(g_0) \asymp_{\mathcal{Q}} \pi(g_1)$, **D6** in $\mathcal{Q}$ gives $\pi([g_0, g_1]) = [\pi(g_0), \pi(g_1)] \prec_{\mathcal{Q}} \pi(g_0)$, whence $[g_0, g_1] \prec g_0$. Suppose now that $g_0 \in \mathcal{N}$. So likewise $g_1 \in \mathcal{N}$, and thus we must have $g_0 \asymp_{\mathcal{N}} g_1$. So **D6** in $\mathcal{N}$ yields $[g_0, g_1] \prec_{\mathcal{N}} g_0$, i.e. $[g_0, g_1] \prec g_0$. This shows that **D6** holds in $\mathcal{G}$, and concludes the proof. □

**Corollary 2.21.** *Let* $(\mathcal{G}, \cdot, 1, \preccurlyeq)$ *and* $(\mathcal{H}, \cdot, 1, \preccurlyeq)$ *be c-valued groups. Let* $\varphi: \mathcal{G} \longrightarrow \mathrm{Aut}(\mathcal{H})$ *be a group homomorphism. Suppose that the following condition is satisfied:*

$$\forall g \in \mathcal{G}^{\neq}, \forall h \in \mathcal{H}^{\neq}, \forall f \in \mathcal{H}, \varphi(g)(h) \neq f\,h\,f^{-1}. \tag{2.5}$$

*Then there is a unique c-dominance relation on the semi-direct product* $\mathcal{H} \rtimes_{\varphi} \mathcal{G}$ *for which the inclusion* $\mathcal{H} \longrightarrow \mathcal{H} \rtimes_{\varphi} \mathcal{G}$ *is an embedding and the projection* $\mathcal{H} \rtimes_{\varphi} \mathcal{G} \twoheadrightarrow \mathcal{G}$ *is a homomorphism. Moreover* $\mathrm{val}(\mathcal{H} \rtimes_{\varphi} \mathcal{G}) = \mathrm{val}(\mathcal{H}) \amalg \mathrm{val}(\mathcal{G})$*.*

**Example 2.22.** Consider the c-valued group $C(\!(t)\!)^{\equiv 1}$ of Example 2.11. We have an action $\varphi : C^{\times} \longrightarrow \mathrm{Aut}(C(\!(t)\!)^{\equiv 1})$ of $C^{\times}$ on $C(\!(t)\!)^{\equiv 1}$ by conjugation $c \mapsto (f \mapsto c\, f(c^{-1}\, t))$. The corresponding semidirect product is naturally isomorphic to the group of $C$-linear automorphisms of the ordered field $C(\!(t)\!)$ (see [36]). The action does not satisfy (2.5) since for instance, for $f = t + t^3 + t^5 + \cdots \in C(\!(t)\!)^{\equiv 1}$, we have $f = \varphi(-1)(f)$. However, given a torsion-free subgroup $G \subseteq C^{\times}$, the condition (2.5) holds for the restriction of $\varphi$ to $G$. Indeed, for a non-trivial $f = t + f_{v(f)}\, t^{v(f)} + \delta$ where $v(\delta) > v(f)$ and $c \in G^{\neq}$, we have $\varphi(c)(f) = t + f_{v(f)}\, c^{1 - v(f)}\, t^{v(f)} + \varepsilon$ for an $\varepsilon$ with $v(\varepsilon) > v(f)$. On the other hand $g \circ f \circ g^{\mathrm{inv}} = t + f_{v(f)}\, t^{v(f)} + \cdots$ for all $g \in C(\!(t)\!)^{\equiv 1}$ by Example 2.11(iii). Since $v(f) > 1$ (because $f$ is non-trivial) and $G$ is torsion-free, it follows that $c^{1 - v(f)} \neq 1$, so $\varphi(c)(f) \neq g \circ f \circ g^{\mathrm{inv}}$.

## 2.5 Residues and scaling elements

In this subsection, we fix a valued group $(\mathcal{G}, \cdot, 1, \preccurlyeq)$ and a ring $A$.

**Proposition 2.23.** *For all* $f, g \in \mathcal{G}^{\neq}$*, the following are equivalent:*

$$f\,g^{-1} \prec f; \qquad f\,g^{-1} \prec g; \qquad g^{-1}\,f \prec f; \qquad g^{-1}\,f \prec g.$$

**Proof.** Suppose that $f g^{-1} \prec f$. So $g^{-1} = f^{-1}(f g^{-1}) \asymp f$ by Proposition 2.4, whence $g \asymp g^{-1} \asymp f$ by **D4**. We deduce that $f g^{-1} \prec g$. In view of **D4**, switching $f$ with $g^{-1}$, we obtain the other implications. $\square$

For $f, g \in \mathcal{G}$, we write $f \sim g$ if one of the statements in Proposition 2.23 holds.

**Corollary 2.24.** *For all $f, g \in \mathcal{G}$, we have $f \sim g \Longleftrightarrow f^{-1} \sim g^{-1}$.*

**Lemma 2.25.** *The relation $\sim$ is an equivalence relation on $\mathcal{G}^{\neq}$ which is thinner than $\asymp$.*

**Proof.** Let $f, g, h \in \mathcal{G}^{\neq}$. We have seen in the proof of Proposition 2.23 that $f \sim g \Longrightarrow f \asymp g$, so $\sim$ is thinner than $\asymp$. The symmetry of $\sim$ follows from Proposition 2.23. We have $f \sim f$ since $1 = f f^{-1} \prec f$ by **D1**, so $\sim$ is reflexive. Suppose that $f \sim g$ and $g \sim h$. We have $f g^{-1} \prec g$ and $g h^{-1} \prec g$, whence $f h^{-1} = (f g^{-1})(g h^{-1}) \prec g$ by **D2**. But $g \asymp f$ so $f h^{-1} \prec f$, i.e. $f \sim h$. Thus $\sim$ is transitive. $\square$

Since the relation $\sim$ is definable in $\mathcal{L}_{\mathrm{vg}}$, the axiom **D3** entails that:

**Lemma 2.26.** *For all $f, g, h \in \mathcal{G}$, we have $f \sim g \Longleftrightarrow h f h^{-1} \sim h g h^{-1}$.*

**Definition 2.27.** *Let $(\mathcal{G}, \cdot, 1, \preccurlyeq)$ be a valued group and let $f, g \in \mathcal{G}^{\neq}$ with $f \succ g$. A* **commutator bound** *of $(f, g)$ is an element $\flat \in \mathrm{val}(\mathcal{G})$ such that for all $f_0, g_0 \in \mathcal{G}$ with $f_0 \sim f$ and $g_0 \sim g$, we have $\mathrm{val}([f_0, g_0]) \geqslant \flat$. We say that $\mathcal{G}$ is* **c-bounded** *if all ordered pairs $(f, g) \in \mathcal{G}^{\neq} \times \mathcal{G}^{\neq}$ with $f \succ g$ have a commutator bound.*

**Remark 2.28.** Note that c-boundedness is preserved under coarsenings.

It is immediate that:

**Lemma 2.29.** *Any c-bounded valued group satisfies* **D5**.

For $f \in \mathcal{G}^{\neq}$, we denote by $\mathrm{res}(f)$ the equivalence class of $f$ for $\sim$, called the *residue* of $f$. We also set $\mathrm{res}(1) := \{1\}$. Let $\rho \in \mathrm{val}(\mathcal{G})$. We define a group $\mathcal{C}_\rho$ as follows. Set

$$\mathcal{C}_\rho := \{\mathrm{res}(f) : f = 1 \vee \mathrm{val}(f) = \rho\}.$$

For $\mathrm{res}(f), \mathrm{res}(g) \in \mathcal{C}_\rho$, set

$$\begin{aligned} \mathrm{res}(f) + \mathrm{res}(g) &:= \mathrm{res}(f g) \quad \text{if } f g \asymp f \\ \mathrm{res}(f) + \mathrm{res}(g) &:= \mathrm{res}(1) \quad \text{if } f g \prec f. \end{aligned}$$

In other words $(\mathcal{C}_\rho, +, \mathrm{res}(1))$ is isomorphic to the quotient group $\mathcal{G}_{\leqslant \rho} / \mathcal{G}_{< \rho}$ of Proposition 2.19. In particular, it is a group. We use additive denotations for $\mathcal{C}_\rho$ because will be Abelian in most applications:

**Proposition 2.30.** *If $\mathcal{G}$ satisfies* **D6***, then $(\mathcal{C}_\rho, +, \mathrm{res}(1))$ is Abelian.*

**Proof.** Let $\mathrm{res}(f), \mathrm{res}(g) \in \mathcal{C}_\rho$. If $f g \prec f$, then $g f = f^{-1}(f g) f \prec f^{-1} f f = f$ by **D3**, so $\mathrm{res}(f) + \mathrm{res}(g) = 0 = \mathrm{res}(g) + \mathrm{res}(f)$. If $f g \asymp f$, then $g f = f g [g, f]$ where $[g, f] \prec f \asymp f g$ by **D6**. Therefore $g f \sim f g$, i.e. $\mathrm{res}(g) + \mathrm{res}(f) = \mathrm{res}(f) + \mathrm{res}(g)$. $\square$

It is easy to see by induction on $n \in \mathbb{N}$ that if $\mathcal{G}$ is c-valued group with $\mathrm{val}(\mathcal{G}) = n$, then it is an iterated extension of Abelian (trivially) c-valued groups. More precisely, writing $\mathrm{val}(\mathcal{G}) = \{\rho_1, \ldots, \rho_n\}$ where $\rho_1 < \cdots < \rho_n$, we see that $\mathcal{G}$ sits in the diagram

$$\begin{array}{ccccccccc} 0 & \hookrightarrow & \mathcal{C}_{\rho_1} \simeq \mathcal{G}_{<\rho_2} & \hookrightarrow & \mathcal{G}_{\leqslant \rho_2} & \hookrightarrow & \cdots & \hookrightarrow & \mathcal{G}_{\leqslant \rho_n} = \mathcal{G} \\ & & & & \downarrow & & & & \downarrow \\ & & & & \mathcal{C}_{\rho_2} & & & & \mathcal{C}_{\rho_n} \\ & & & & \downarrow & & & & \downarrow \\ & & & & 0 & & & & 0 \end{array}$$

of valued group homomorphisms where the upper arrows are monomorphisms and the ⌝-shaped sequences of three groups are exact.

Given a morphism (resp. an embedding) of valued groups $\Phi : \mathcal{H} \longrightarrow \mathcal{G}$ and a $\rho = \mathrm{val}(h) \in \mathrm{val}(\mathcal{H})$, we have a group homomorphism (resp. embedding)

$$\begin{array}{rcl} \mathrm{res}_{\Phi} : \mathcal{C}_{\rho} & \longrightarrow & \mathcal{C}_{\mathrm{val}(\Phi(h))} \\ \mathrm{res}(f) & \longmapsto & \mathrm{res}(\Phi(f)). \end{array}$$

**Proposition 2.31.** *Let $n > 0$ and $f_1, \ldots, f_n \in \mathcal{G}$. Writing*

$$J = \{ i \in \{1, \ldots, n\} : \mathrm{val}(f_i) = \max(\mathrm{val}(f_1), \ldots, \mathrm{val}(f_n)) \},$$

*we have* $\mathrm{res}(f_1 \cdots f_n) = \sum_{j \in J} \mathrm{res}(f_j)$ *if* $\sum_{j \in J} \mathrm{res}(f_j) \neq 0$ *and* $f_1 \cdots f_n \prec f_j$ *for all* $j \in J$ *otherwise.*

**Proof.** Let $n, f_1, \ldots, f_n, J$ be as in the statement. The result is trivial if $f_1 = \cdots = f_n = 1$ so we may assume that at least one of the $f_i$'s is non-trivial. We prove the result by induction on $n$. If $n = 0$ or $n = 1$, this is immediate. Suppose that $n \geqslant 2$ and that the statement holds for $n - 1$. Write

$$\begin{array}{rcl} J' & := & \{ i \in \{1, \ldots, n-1\} : \mathrm{val}(f_i) = \max(\mathrm{val}(f_1), \ldots, \mathrm{val}(f_{n-1})) \} \\ J & := & \{ i \in \{1, \ldots, n\} : \mathrm{val}(f_i) = \max(\mathrm{val}(f_1), \ldots, \mathrm{val}(f_n)) \}, \end{array}$$

and fix a $j \in J'$. If $f_n \asymp f_j$ and $\sum_{i \in J'} \mathrm{res}(f_i) = 0$, then since $f_j \succ 1$ we have $\mathrm{res}(f_n) = \sum_{i \in J} \mathrm{res}(f_i) \neq 0$, and $f_1 \cdots f_{n-1} \prec f_j \asymp f_n$ by the induction hypothesis. It follows that $f_1 \cdots f_n \sim f_n$, i.e. $\mathrm{res}(f_1 \cdots f_n) = \mathrm{res}(f_n) = \sum_{i \in J} \mathrm{res}(f_i)$. If $f_n \asymp f_j$ and $\sum_{i \in J'} \mathrm{res}(f_i) = -\mathrm{res}(f_n)$, then the induction hypothesis gives $f_1 \cdots f_{n-1} \sim f_n^{-1}$, so $f_1 \cdots f_n \prec f_n \asymp f_j$. Finally if $f_n \asymp f_j$ and

$$c := \sum_{i \in J'} \mathrm{res}(f_i) \notin \{0, -\mathrm{res}(f_n)\},$$

then the induction hypothesis gives $\mathrm{res}(f_1 \cdots f_{n-1}) = c$, whence in particular $f_1 \cdots f_{n-1} \asymp f_n$. We deduce since $c + \mathrm{res}(f_n) \neq 0$ that $f_1 \cdots f_n \not\prec f_n$, so $\mathrm{res}(f_1 \cdots f_{n-1} f_n) = \mathrm{res}(f_1 \cdots f_{n-1}) + \mathrm{res}(f_n) = \sum_{i \in J} \mathrm{res}(f_i)$ as claimed.

If $f_n \succ f_j$, then $J = \{n\}$ and $f_1 \cdots f_{n-1} \preccurlyeq f_j$ by the induction hypothesis. Thus $f_1 \cdots f_n \sim f_n$, so $\mathrm{res}(f_1 \cdots f_n) = \mathrm{res}(f_n) = \sum_{i \in J} \mathrm{res}(f_i)$. Lastly, if $f_n \prec f_j$, then $J = J'$, and either $\sum_{i \in J'} \mathrm{res}(f_i) = 0$ and $f_1 \cdots f_{n-1} \prec f_j$ and $f_1 \cdots f_n \prec f_j$ as well, or $\sum_{i \in J'} \mathrm{res}(f_i) \neq 0$, and then $f_1 \cdots f_{n-1} \asymp f_j$ by the induction hypothesis. In that case $f_1 \cdots f_n \sim f_1 \cdots f_{n-1}$ so $\mathrm{res}(f_1 \cdots f_n) = \sum_{i \in J} \mathrm{res}(f_i)$. We conclude by induction. □

**Definition 2.32.** *An $\mathscr{s} \in \mathcal{G}^{\neq}$ is said* **scaling** *if $\mathcal{C}(\mathscr{s})$ is Abelian and for all $f \asymp \mathscr{s}$, there is an $\mathscr{s}_0 \in \mathcal{C}(\mathscr{s})$ with $\mathscr{s}_0 \sim f$.*

**Proposition 2.33.** *Suppose that $\mathcal{G}$ is c-valued. For $\mathscr{s} \in \mathcal{G}^{\neq}$, the following assertions are equivalent:*

*a) $\mathscr{s}$ is scaling;*

b) *the function* $\operatorname{res}:\mathcal{C}(\jmath)\longrightarrow\mathcal{C}_{\operatorname{val}(\jmath)}$ *is a group isomorphism;*

c) *the group* $\mathcal{C}(\jmath)$ *is a complement of* $\mathcal{G}_{<\operatorname{val}(\jmath)}$ *in* $\mathcal{G}_{\leqslant\operatorname{val}(\jmath)}$.

**Proof.** That a and b are equivalent is clear. For any $\jmath\in\mathcal{G}^{\neq}$, we have $\mathcal{C}(\jmath)\cap\mathcal{G}_{<\operatorname{val}(\jmath)}=\{1\}$ by **D5**. Thus $\mathcal{C}(\jmath)$ is a complement of $\mathcal{G}_{<\operatorname{val}(\jmath)}$ in $\mathcal{G}_{\leqslant\operatorname{val}(\jmath)}$ if and only if in $\mathcal{G}_{\leqslant\operatorname{val}(\jmath)}=\bigcup_{f\in\mathcal{C}(\jmath)}f\mathcal{G}_{<\operatorname{val}(\jmath)}=\bigcup_{f\in\mathcal{C}(\jmath)}\operatorname{res}(f)=\operatorname{res}(\mathcal{C}(\jmath))$. So c and b are equivalent. □

**Lemma 2.34.** *Let* $C$ *be a field of characteristic zero and consider the the c-valued group* $C(\!(t)\!)^{\equiv 1}$ *of Example 2.11. Any non-trivial element* $f=t+\sum_{n>1}f_n t^n$ *of* $C(\!(t)\!)^{\equiv 1}$ *is scaling, and* $\mathcal{C}(f)$ *is isomorphic to the underlying additive group of* $C$.

**Proof.** For all $\mu\in C$, there is [9, Proposition 2.16] a series $f^{[\mu]}\in\mathcal{C}(f)$ with $f^{[\mu]}=t+\mu f_{\operatorname{val}(f)}t^{\operatorname{val}(f)}+\cdots$. So given $g=t+\gamma t^{\operatorname{val}(f)}+\cdots$ with $\operatorname{val}(g)=\operatorname{val}(f)$ (i.e. $\gamma\neq 0$), we have $f^{[\gamma/\lambda]}\sim g$. So $f$ is scaling and $\mu\mapsto f^{[\mu]}$ is a group embedding by [9, Proposition 2.16]. Since $C\longrightarrow\mathcal{C}_{\operatorname{val}(f)}\,;\mu\mapsto\operatorname{res}(f^{[\mu]})$ is surjective, this embedding is surjective. □

# 3 Valued exponential groups

Our rings are always associative and unital. Throughout this section, we fix a ring $A$.

## 3.1 Exponential groups

As per Miasnikov-Remeslennikov [26, 27], an $A$-*group* is a group $\mathcal{G}$ together with a map $A\times\mathcal{G}\longrightarrow\mathcal{G}\,;(a,g)\mapsto g^a$ satisfying the following axioms for all $g,h\in\mathcal{G}$ and $a,b\in A$:

1. $g^0=1$, $g^1=g$, $1^a=1$.
2. $g^{a+b}=g^a g^b$, $(g^a)^b=g^{ab}$.
3. $(h\,g\,h^{-1})^a=h\,g^a\,h^{-1}$.
4. $[g,h]=1\rightarrow(g\,h)^a=g^a\,h^a$.

We see $A$-groups as one-sorted structures $(\mathcal{G},\cdot,1,(g\mapsto g^a)_{a\in A})$, and we denote by $\mathcal{L}_{A\mathrm{g}}$ the first-order language of $A$-groups. An $A$-*subgroup* of $\mathcal{G}$ is a subgroup which is closed under the $a$-power maps $g\mapsto g^a, a\in A$, i.e. an $\mathcal{L}_{A\mathrm{g}}$-substructure of $\mathcal{G}$. An $A$-*ideal* of $\mathcal{G}$ is a normal $A$-subgroup $\mathcal{H}$ such that for all $f,g\in\mathcal{G}$ and $a\in A$, we have $f^{-1}g^{-1}fg\in\mathcal{H}\Longrightarrow f^{-a}g^{-a}f^a g^a\in\mathcal{H}$. Equivalently [26, Proposition 5], it is the kernel of an $A$-*homomorphism* (i.e. $\mathcal{L}_{A\mathrm{g}}$-homomorphism) of $A$-groups. The $A$-closure $\operatorname{cl}_A(X)$ of a subset $X\subseteq\mathcal{G}$ is the smallest $A$-ideal of $\mathcal{G}$ containing $X$. It is obtained as an increasing union $\operatorname{cl}_A(X)=\bigcup_{n\in\mathbb{N}}X_n$ where $X_0=X$ and each $X_{n+1}$ is generated by conjugates and $a$-powers of elements of $X_n$, as well as by elements $f^{-a}g^{-a}f^a g^a$ for $f^{-1}g^{-1}fg\in X_n$, where $a$ ranges in $A$.

The derived $A$-ideal of an $A$-group $\mathcal{G}$ is defined as the $A$-closure $\operatorname{cl}_A([\mathcal{G},\mathcal{G}])$ of the derived subgroup $[\mathcal{G},\mathcal{G}]$ of $\mathcal{G}$, and denoted $\mathcal{G}'$. The (transfinite) *derived* $A$-*series* of $\mathcal{G}$ is the sequence $(\mathcal{G}^{(\gamma)})_{\gamma\in\mathbf{On}}$ of $A$-ideals of $\mathcal{G}$ where $\mathcal{G}^{(0)}=\mathcal{G}$, $\mathcal{G}^{(\gamma+1)}=(\mathcal{G}^{(\gamma)})'$ for all $\gamma\in\mathbf{On}$ and $\mathcal{G}^{(\lambda)}=\bigcap_{\gamma<\lambda}\mathcal{G}^{(\gamma)}$ for non-zero limit ordinals $\lambda$. The $A$-group $\mathcal{G}$ is said $A$-*hypoabelian* if $\mathcal{G}^{(\gamma)}=\{1\}$ for some $\gamma\in\mathbf{On}$. If $\gamma\leqslant\omega$, then we say that $\mathcal{G}$ is residually $A$-solvable.

Likewise, the *lower central* $A$-*series* $(\mathcal{G}_\gamma)_{\gamma\in\mathbf{On}}$ is defined by $\mathcal{G}_0:=\mathcal{G}$, $\mathcal{G}_{\gamma+1}:=\operatorname{cl}_A([\mathcal{G},\mathcal{G}_\gamma])$ and $\mathcal{G}_\lambda:=\bigcap_{\gamma<\lambda}\mathcal{G}_\gamma$ if $\lambda$ is a non-zero limit. The $A$-group $\mathcal{G}$ is said $A$-*hypocentral* if $\mathcal{G}_\gamma=\{1\}$ for some $\gamma\in\mathbf{On}$. If $\gamma\leqslant\omega$, then we say that $\mathcal{G}$ is residually $A$-nilpotent.

Lastly, we say that $\mathcal{G}$ is $A$-*torsion-free* if we have $g^a\neq 1$ for all $a\in A\setminus\{0\}$ and $g\in\mathcal{G}^{\neq}$.

**Remark 3.1.** Any group $(\mathcal{G},\cdot,1)$ has a unique structure of $\mathbb{Z}$-group, given by $(n,g)\mapsto g^n$. All $\mathbb{Z}$-notions coincide with the usual purely group theoretic ones.

**Remark 3.2.** That the elements of the derived and lower central $A$-series are $A$-ideals of $\mathcal{G}$ itself follows by induction as in the pure group case, using the increasing union presentation of $\mathrm{cl}_A(X)$ for $X\subseteq\mathcal{G}$.

## 3.2 Valuations on exponential groups

We write $\mathcal{L}_{A\mathrm{vg}}$ for the first-order language of $A$-groups expanded with a symbol of binary relation $\preccurlyeq$.

**Definition 3.3.** *Given an $A$-group $(\mathcal{G},\cdot,1,(g\mapsto g^a)_{a\in A})$ and a dominance relation $\preccurlyeq$ on $\mathcal{G}$, consider the following axioms in $\mathcal{L}_{A\mathrm{vg}}$:*

**D7$^A$.** *For all $f\in\mathcal{G}$ and $a\in A$, we have $f^a\neq 1\Longrightarrow f^a\asymp f$.*

**D8$^A$.** *For all $f,\varepsilon\in\mathcal{G}^{\neq}$ and $a\in A$, we have $\varepsilon\prec f\Longrightarrow(f\varepsilon)^a f^{-a}\prec f$.*

**D9$^A$.** *For all $f,\varepsilon\in\mathcal{G}^{\neq}$ and $a\in A$, we have $\varepsilon\prec f\Longrightarrow(f\varepsilon)^a f^{-a}\preccurlyeq\varepsilon$.*

*We say that $(\mathcal{G},\cdot,1,(g\mapsto g^a)_{a\in A},\preccurlyeq)$ is* ***A*-valued** *if* **D7$^A$** *and* **D8$^A$** *hold. We say that it is* **uniformly *A*-valued** *if* **D7$^A$** *and* **D9$^A$** *hold.*

We call the corresponding valuation an $A$-valuation. Note that **D7$^A$** is weaker than **D5** for a valued group. Furthermore, if $A=\mathbb{Z}$, then **D8$^A$** follows by induction from Lemma 2.5 and Proposition 2.4. Therefore:

**Proposition 3.4.** *Any c-valued group is a $\mathbb{Z}$-valued $\mathbb{Z}$-group.*

## 3.3 Examples of valued exponential groups

**Definition 3.5.** *Let $\lambda$ be an ordinal and let $(\mathcal{G}_{[\gamma]})_{\gamma\leqslant\lambda}$ be decreasing sequence of $A$-ideals of $\mathcal{G}_{[0]}$. We say that $(\mathcal{G}_{[\gamma]})_{\gamma\leqslant\beta}$ is an* **A-series** *if*

- *$\mathcal{G}_{[\beta]}=\{1\}$;*
- *$\mathcal{G}_{[\lambda]}=\bigcap_{\gamma<\lambda}\mathcal{G}_{[\gamma]}$ for all non-zero limit ordinals $\lambda\leqslant\beta$;*
- *each quotient $\mathcal{G}_{[\gamma]}/\mathcal{G}_{[\gamma+1]}$ for $\gamma<\beta$ is an (Abelian) $A$-module without $A$-torsion.*

Given an $A$-series $(\mathcal{G}_{[\gamma]})_{\gamma\leqslant\beta}$, we define an associated valuation $v$ on $\mathcal{G}_{[0]}$ with value set $v(\mathcal{G}_{[0]})=(\beta,>)$, by setting $v(g)$ to be the maximal ordinal with $g\in\mathcal{G}_{[v(g)]}$ for all $g\in\mathcal{G}_{[0]}^{\neq}$. We call the corresponding dominance relation the associated dominance relation.

**Proposition 3.6.** *Let $(\mathcal{G}_{[\gamma]})_{\gamma\leqslant\beta}$ be an $A$-series. Then $\mathcal{G}_{[0]}$ is uniformly $A$-valued for the associated dominance relation. Moreover* **D6** *holds.*

**Proof.** The axiom **D1** trivially holds. Each term in the $A$-series is a normal subgroup, so the axioms **D2**, **D3** and **D4**. Since each quotient $A$-module in the $A$-series has no $A$-torsion, **D7$^A$** holds. Since each term in the $A$-series is an $A$-ideal of $\mathcal{G}_{[0]}$, **D9$^A$** also holds. Lastly **D6** holds because quotients in the series are Abelian. □

**Example 3.7.** If $\mathcal{G}$ is an $A$-hypoabelian (resp. $A$-hypocentral) $A$-group, then the derived (resp. lower central) $A$-series of $\mathcal{G}$ is an $A$-series. We call the associated valuation the *derived valuation* (resp. *lower central $A$-valuation*) on $\mathcal{G}$.

**Remark 3.8.** The axiom **V5** may not hold for the associated valuation. For instance, let $\mathcal{H},\mathcal{G}$ be two non-trivial solvable $\mathbb{Z}$-groups and assume that $\mathcal{G}$ is non-Abelian. It is easy to see that for the derived valuation $v$ on $\mathcal{H}\times\mathcal{G}$, for all all $(h,g)\in\mathcal{H}\times\mathcal{G}$, we have $v(h,g)=\min(v(h),v(g))$. The conditions on $\mathcal{H},\mathcal{G}$ imply that there is an $(h,g)\in\mathcal{H}^{\neq}\times\mathcal{G}^{\neq}$ with $v(h)\neq v(g)$. But then $(h,1)$ and $(1,g)$ commute and have distinct natural valuations, so **V5** does not hold.

We recall [27, Definition 3] that a *free $A$-group* is an $A$-group $\mathcal{F}$ such that there is a subset $X\subseteq\mathcal{F}$ such that for all $A$-groups $\mathcal{G}$, all maps $\varphi: X\longrightarrow\mathcal{G}$ extend uniquely into $A$-group homomorphisms $\mathcal{F}\longrightarrow\mathcal{G}$. We say that $\mathcal{F}$ is a free $A$-group over $X$. For any set $X$, there is a unique free $A$-group over $X$ up to $A$-group isomorphism.

**Proposition 3.9.** *Suppose that $A$ is a commutative domain of characteristic zero. Let $\mathcal{F}$ be a free $A$-group. Then the derived and lower central valuations on $\mathcal{F}$ are c-valuations.*

**Proof.** By [21, Corollary 5.9], the group $\mathcal{F}$ is residually $A$-torsion-free nilpotent, so the lower central and derived valuations are defined. We claim that the lower central valuation $v$ on $\mathcal{F}$ satisfies **V5**. Indeed let $f,g\in\mathcal{F}^{\neq}$ with $g\in\mathcal{C}(f)$. By [27, Proposition 8.16], we have $\mathcal{C}(f)=s^A$ for a certain $s\in\mathcal{C}(f)$. Let $a,b\in A\setminus\{0\}$ with $(f,g)=(s^a,s^b)$. It suffices to show that $v(f)=v(s)$ and that $v(g)=v(s)$. Set $n:=v(s)$. The $A$-group $\mathcal{F}_n$ is an $A$-subgroup of $\mathcal{F}$, hence it is a free $A$-group as a consequence of [27, Theorem 7.8]. Assume for contradiction that $s^a\in\mathcal{F}_{n+1}$. Then $(s\mathcal{F}_{n+1})^a=1$ in the Abelian $A$-module $\mathcal{F}_n/\mathcal{F}_{n+1}$. But $\mathcal{F}_n/\mathcal{F}_{n+1}$ is a free Abelian $A$-module. Since $A$ has no zero-divisors, it follows that $s\in\mathcal{F}_{n+1}$: a contradiction. So $f\notin\mathcal{F}_{n+1}$ and likewise $g\notin\mathcal{F}_{n+1}$, i.e. $v(f)=n=v(g)$. Thus **V5** holds, so $(\mathcal{F},\cdot,1,v)$ is c-valued. The same arguments apply for the derived valuation. □

## 3.4 Residues as exponential groups

**Proposition 3.10.** *Let $\mathcal{G}$ be an $A$-group. Let $\mathcal{N}\trianglelefteq\mathcal{G}$ be an $A$-ideal which is $\preccurlyeq$-initial. Then $\mathcal{G}/\mathcal{N}$ satisfies* **D7**$^A$ *(resp.* **D8**$^A$*, resp.* **D9**$^A$*) if $\mathcal{G}$ does.*

**Proof.** Suppose that $\mathcal{G}$ satisfies **D7**$^A$. Let $f\in\mathcal{G}$ and $a\in A$ with $(f\mathcal{N})^a=f^a\mathcal{N}\neq\mathcal{N}$. We have $f,f^a\notin\mathcal{N}$, so $f^a\neq 1$. We deduce by **D7**$^A$ that $f\asymp f^a$, whence also $f\mathcal{N}\asymp(f\mathcal{N})^a$. So **D7**$^A$ holds in $\mathcal{G}/\mathcal{N}$. Suppose now that $\mathcal{G}$ satisfies **D8**$^A$ (resp. **D9**$^A$). Let $f,\varepsilon\in\mathcal{G}\setminus\mathcal{N}$ with $f\mathcal{N}\succ\varepsilon\mathcal{N}$, and let $a\in A$. We have $f\succ\varepsilon$, so **D8**$^A$ (resp. **D9**$^A$) in $\mathcal{G}$ gives $(f\varepsilon)^a f^{-a}\prec f$ (resp. $(f\varepsilon)^a f^{-a}\preccurlyeq\varepsilon$), whence also

$$(f\mathcal{N}\varepsilon\mathcal{N})^a(f\mathcal{N})^{-a}=(f\varepsilon)^a f^{-a}\mathcal{N}\prec f\mathcal{N}\qquad(\text{resp. }(f\mathcal{N}\varepsilon\mathcal{N})^a(f\mathcal{N})^{-a}\preccurlyeq\varepsilon\mathcal{N})$$

So **D8**$^A$ (resp. **D9**$^A$) holds in $\mathcal{G}/\mathcal{N}$. □

**Remark 3.11.** Suppose that $\mathcal{G}$ is an $A$-valued $A$-group. Then by **D7**$^A$, the group $\mathcal{G}_{\leqslant\rho}$ is $A$-valued for the induced structure, and by **D8**$^A$, the subgroup $\mathcal{G}_{<\rho}$ is an $A$-ideal of $\mathcal{G}_{\leqslant\rho}$. Thus $\mathcal{C}_\rho$ is an $A$-group for the operation $A\times\mathcal{C}_\rho\longrightarrow\mathcal{C}_\rho$; $(a,\operatorname{res}(f))\mapsto a\operatorname{res}(f):=\operatorname{res}(f^a)$.

**Lemma 3.12.** *Suppose that $\mathcal{G}$ is an $A$-valued $A$-group. Let $a\in A$ and $\operatorname{res}(f)\in\mathcal{C}_\rho$, we have $a\operatorname{res}(f)=0$ if and only if $f^a=1$.*

**Proof.** If $f^a=1$, then $a\operatorname{res}(f)=\operatorname{res}(f^a)=0$. Conversely, if $a\operatorname{res}(f)=0$, then $f^a\prec f$, which by virtue of **D7**$^A$ implies that $f^a=1$. □

Therefore $\mathcal{G}$ is $A$-torsion-free if and only if all $\mathcal{C}_\rho$ are $A$-torsion-free.

# 4 Ordered valued groups

## 4.1 Ordered valued groups

An *ordered group* is a group $(\mathcal{G}, \cdot, 1)$ together with a total ordering $<$ on $\mathcal{G}$ such that

$$\forall f, g, h \in \mathcal{G}, g > h \Longrightarrow (f g > f h \wedge g f > h f). \tag{4.1}$$

We write $\leqslant$ for the relation $f \leqslant g \Longleftrightarrow (f < g \vee f = g)$. We have a first-order language of ordered groups $\mathcal{L}_{\mathrm{og}} := \langle \cdot, 1, \leqslant, \mathrm{Inv} \rangle$. Given an ordered group $\mathcal{G}$, we write

$$\mathcal{G}^{>} := \{f \in \mathcal{G} : f > 1\} \qquad \text{and} \qquad \mathcal{G}^{\geqslant} := \{f \in \mathcal{G} : f \geqslant 1\}.$$

**Remark 4.1.** All ordered groups are torsion-free [25]. They also satisfy [29, Lemma 1.1] the consequence of Lemma 2.14.

**Definition 4.2.** *An* **ordered (c-)valued group** *is an ordered group* $(\mathcal{G}, \cdot, 1, <)$ *together with a (c-)dominance relation* $\preccurlyeq$ *such that each set* $\mathrm{val}(g) \cap \mathcal{G}^{>}$, $g \in \mathcal{G}^{>}$ *is convex.*

Equivalently, an ordered (c-)valued group is an ordered group $\mathcal{G}$ together with a (c-)valuation that is nondecreasing on $\mathcal{G}^{>}$, or again equivalently an ordered group with a (c-)dominance relation whose intersection with $\mathcal{G}^{>}$ is coarser than $\leqslant$. We write $\mathcal{L}_{\mathrm{ovg}}$ for the expansion of $\mathcal{L}_{\mathrm{og}}$ with a binary relation symbol interpreted as the dominance relation in ordered valued groups. The axioms of ordered (c-)valued groups are universal, so:

**Proposition 4.3.** *A subgroup of an ordered (c-)valued group is an ordered (c-)valued group for the induced structure.*

**Example 4.4.** Any ordered Abelian group is an ordered c-valued group for its unique c-dominance relation (see Remark 2.7).

**Remark 4.5.** If $(\mathcal{G}, \cdot, 1, <, \preccurlyeq)$ is an ordered valued group, then the inverse group $(\check{\mathcal{G}}, \check{\cdot}, 1, <, \preccurlyeq)$ is also an ordered valued group.

Let us give an important example of ordered valued groups. Let $(K, +, \cdot, 0, 1, \mathcal{O}, \partial, \circ, x)$ be an H-field with composition and inversion as per [8, Section 4.1], let $\mathcal{o}$ be the maximal ideal of the valuation ring $\mathcal{O}$ and let $C = \mathrm{Ker}(\partial)$ be the field of constants. We write $a^{\dagger} := \frac{\partial(a)}{a}$ for each $f \in K^{\times}$. By definition $K$ satisfies the Taylor axiom scheme which is the following set of sentences where $n$ ranges in $\mathbb{N}$:

**HFC5.** $\forall f, g, \delta \in K, g > C \wedge \delta \in g\, \mathcal{o} \Rightarrow f - \sum_{i \leqslant n} \frac{f^{(i)} \circ g}{i!} \delta^i \in ((f^{(n)} \circ g)\, \delta^n)\, \mathcal{o}$.

Let us show that $K$ induces an ordered valued group.

**Lemma 4.6.** *Let* $a \in K^{>C}$ *and* $\delta \in x\, \mathcal{o}$. *We have* $a \circ (x + \delta) - a \in a\, \mathcal{o} \Longleftrightarrow a^{\dagger} \delta \in \mathcal{o}$.

**Proof.** If $b^{\dagger} \delta \in \mathcal{o}$, then **HFC5** gives $a \circ (x + \delta) - a + a' \delta \in (a' \delta)\, \mathcal{o}$ where $a' \delta \in a\, \mathcal{o}$, hence the result. Suppose now that $a^{\dagger} \delta \notin \mathcal{o}$. We also assume that $\delta > 0$, the case when $\delta < 0$ being symmetric. Let $\mathfrak{m} \in K$ with $\mathfrak{m} > 0$ and $a^{\dagger} \mathfrak{m} \prec 1$. Then $\mathfrak{m} \prec \delta$ so [8, Section 4.1, **HFC3**] gives $a \circ (x + \delta) - a > a \circ (x + \mathfrak{m}) - a$. By **HFC5**, this gives $a \circ (x + \delta) - a > a' \mathfrak{m} + \iota$ for a $\iota \in (a' \mathfrak{m})\, \mathcal{o}$. So $a \circ (x + \delta) - a > a' \mathfrak{m}$ pour tout $\mathfrak{m} \in K$ with $\mathfrak{m} > 0$ and $a^{\dagger} \mathfrak{m} \prec 1$. This entails that $a \circ (x + \delta) - a \notin (a' a^{\dagger})\, \mathcal{o} = a\, \mathcal{o}$. $\square$

Consider the set $\mathcal{P} := \{x + \delta \in K : \delta \in x\,\mathcal{o}\}$ of so-called *parabolic* elements in $K$ (see [30]), with the induced ordering $<$ and composition law $\circ$. We can see by **HFC5** and Lemma 4.6 that the inverse $(x+\delta)^{\mathrm{inv}}$ of an $x+\delta \in \mathcal{T}$ satisfies

$$(x+\delta)^{\mathrm{inv}} - x - \delta \in \delta\,\mathcal{o}. \tag{4.2}$$

So $\mathcal{P}$ is a subgroup of $(K^{>C}, \circ, x)$. We endow $\mathcal{P}$ with the total quasi-ordering

$$\forall \varepsilon, \delta \in x\,\mathcal{o}, (x + \varepsilon \preccurlyeq x + \delta \Longleftrightarrow \varepsilon \in \delta\,\mathcal{O}).$$

**Lemma 4.7.** *For all $x+\delta, x+\varepsilon \in \mathcal{P}$, we have*

$$(x+\delta) \circ (x+\varepsilon) \circ (x+\delta)^{\mathrm{inv}} - x + \varepsilon \circ (x+\delta)^{\mathrm{inv}} \in (\varepsilon \circ (x+\delta)^{\mathrm{inv}})\,\mathcal{o}.$$

**Proof.** We have $(x+\delta) \circ (x+\varepsilon) \circ (x+\delta)^{\mathrm{inv}} = (x+\delta) \circ ((x+\delta)^{\mathrm{inv}} + \varepsilon \circ (x+\delta)^{\mathrm{inv}})$. Now $((x+\delta)^{\dagger} \circ (x+\delta)^{\mathrm{inv}})\,(\varepsilon \circ (x+\delta)^{\mathrm{inv}}) \in \mathcal{o}$ if and only if $(x+\delta)^{\dagger}\,\varepsilon \in \mathcal{o}$ by [8, Section 4.1, **HFC1**]. But $(x+\delta)^{\dagger}\,\varepsilon \in x^{\dagger}\,\varepsilon + (x^{\dagger}\,\varepsilon)\,\mathcal{o}$ where $x^{\dagger}\,\varepsilon = \frac{\varepsilon}{x} \in \mathcal{o}$. So we may apply **HFC5** and obtain

$$\begin{aligned}(x+\delta) \circ (x+\varepsilon) \circ (x+\delta)^{\mathrm{inv}} &= (x+\delta) \circ (x+\delta)^{\mathrm{inv}} + ((x+\delta)'\,\varepsilon) \circ (x+\delta)^{\mathrm{inv}} + \iota_0 \\ &= x + \varepsilon \circ (x+\delta)^{\mathrm{inv}} + \iota_1\end{aligned}$$

for $\iota_0, \iota_1 \in \varepsilon \circ (x+\delta)^{\mathrm{inv}}\,\mathcal{o}$ as claimed. □

**Proposition 4.8.** *The structure $(\mathcal{P}, \circ, x, <, \preccurlyeq)$ is an ordered valued group.*

**Proof.** We have $0\,\mathcal{O} = \{0\}$ so **D1** holds. Let $x+\delta, x+\varepsilon \in \mathcal{P}$ and set $f = \max(x+\delta, (x+\delta)^{\mathrm{inv}}, x+\varepsilon, (x+\varepsilon)^{\mathrm{inv}})$. We have $f^{\mathrm{inv}} \circ f^{\mathrm{inv}} \leqslant (x+\delta) \circ (x+\varepsilon) \leqslant f \circ f$ since $\mathcal{P}$ is an ordered group. Now by **HFC5** that $f \circ f = x + 2\,\mu + \iota$ for a $\mu \in \{\delta, \varepsilon, -\delta, -\varepsilon\}$ and a $\iota \in \mu\,\mathcal{o}$. The relation $\preccurlyeq$ is clearly coarser than $<$ on $\mathcal{P}^{>}$. We deduce since $2\,\mu \in \mu\,\mathcal{O}$ that $f \circ f \preccurlyeq x + \delta$ or $f \circ f \preccurlyeq x + \varepsilon$. So **D2** holds. We see with Lemma 4.7 and [8, Section 4.1, **HFC1**] that **D3** holds, whereas **D4** follows from (4.2). □

## 4.2 The growth axiom

**Definition 4.9.** *The* **growth axiom** *is the following $\mathcal{L}_{\mathrm{ovg}}$-sentence:*

**GA.** $\forall f \forall g ((f > 1 \wedge g > 1 \wedge f \succ g) \longrightarrow f\,g\,f^{-1} > g)$.

**Remark 4.10.** Contrary to all axioms introduced so far, **GA** is not preserved under taking the inverse ordered group, or equivalently under reversing the ordering on a given ordered valued group. However, it is preserved under taking subgroups.

**Lemma 4.11.** *Let $(\mathcal{G}, \cdot, 1, <, \preccurlyeq)$ be a an ordered valued group satisfying* **GA**. *For all $f$, $g, h \in \mathcal{G}^{>}$ with $f > h$, $f \not\sim h$ and $h \succ g$, we have $[f^{-1}, g^{-1}] \succcurlyeq [h^{-1}, g^{-1}]$.*

**Proof.** Since $f > h$ and $f \not\sim h$, the element $h^{-1}\,f$ is positive with $h^{-1}\,f \succcurlyeq f \succ g$. Therefore **GA** gives $(h^{-1}\,f)\,g\,(h^{-1}\,f)^{-1} > g$, i.e. $f\,g\,f^{-1} > h\,g\,h^{-1}$, whence $[f^{-1}, g^{-1}] > [h^{-1}, g^{-1}]$. We also have $[h^{-1}, g^{-1}] > 1$ by **GA**. Since $\preccurlyeq$ is coarser than $\leqslant$ on $\mathcal{G}^{\geqslant}$, we deduce that $[f^{-1}, g^{-1}] \succcurlyeq [h^{-1}, g^{-1}]$. □

**Lemma 4.12.** *Let $(\mathcal{G}, \cdot, 1, <, \preccurlyeq)$ be a an ordered valued group satisfying* **GA**. *For all $h$, $g, \varepsilon \in \mathcal{G}^{>}$ with $h \succ g$ and $\varepsilon < g$, we have $[h^{-1}, g^{-1}] \succcurlyeq [h^{-1}, \varepsilon^{-1}]$.*

**Proof.** Note that $\varepsilon \preccurlyeq g$ since $\mathcal{G}$ is an ordered valued group. We have $(\varepsilon^{-1} g)^{-1} < 1$ and $(\varepsilon^{-1} g)^{-1} \prec h$ so **GA** gives $\varepsilon^{-1} g\, h\, (\varepsilon^{-1} g)^{-1} < h$, whence $\varepsilon^{-1} g\, h^{-1} (\varepsilon^{-1} g)^{-1} > h^{-1}$. So $g\, h^{-1} g^{-1} > \varepsilon\, h^{-1} \varepsilon^{-1}$. Multiplying by $h$ on the left, we obtain $[h^{-1}, g^{-1}] \geqslant [h^{-1}, \varepsilon^{-1}]$. We have $[h^{-1}, \varepsilon^{-1}] > 1$ by **GA** so this entails that $[h^{-1}, g^{-1}] \succcurlyeq [h^{-1}, \varepsilon^{-1}]$. □

## 4.3 Ordered quotients

We fix an ordered valued group $(\mathcal{G}, \cdot, 1, <, \preccurlyeq)$. If $\mathcal{N} \trianglelefteq \mathcal{G}$ is a normal convex subgroup of $\mathcal{G}$, then there is a unique ordering on $\mathcal{G}/\mathcal{N}$ for which the quotient map $\mathcal{G} \longrightarrow \mathcal{G}/\mathcal{N}$ is an ordered group homomorphism. It is given by $f\mathcal{N} <_{\mathcal{N}} g\mathcal{N} \Longleftrightarrow f\mathcal{N} < g$ for all $f, g \in \mathcal{N}$. See [25, paragraph 4] for more details.

**Proposition 4.13.** *Let $\mathcal{N} \trianglelefteq \mathcal{G}$ be a normal and $\preccurlyeq$-initial subgroup satisfying the condition $\forall f, g \in \mathcal{G} \setminus \mathcal{N}, [f, g] \in \mathcal{N} \Longrightarrow f \asymp g$. Then $\mathcal{N}$ is convex and $(\mathcal{G}/\mathcal{N}, \cdot, \mathcal{N}, <_{\mathcal{N}})$ is an ordered valued group for the dominance relation of Proposition 2.18. Moreover, it satisfies* **GA** *if $\mathcal{G}$ satisfies* **GA**.

**Proof.** That $\mathcal{N}$ is convex follows from the fact that it is $\preccurlyeq$-initial and that $\mathcal{G}$ is an ordered c-valued group. So $(\mathcal{G}/\mathcal{N}, \cdot, \mathcal{N}, <_{\mathcal{N}})$ is an ordered group. Recall that $(\mathcal{G}/\mathcal{N}, \cdot, \mathcal{N}, \preccurlyeq_{\mathcal{N}})$ is a c-valued group for the c-dominance relation $f\mathcal{N} \preccurlyeq_{\mathcal{N}} g\mathcal{N} \Longleftrightarrow (f \in \mathcal{N} \vee f \preccurlyeq g)$. Given $g \in \mathcal{G}$, it suffices to show that $\operatorname{val}(g\mathcal{N}) \cap (\mathcal{G}/\mathcal{N})^{>}$ is convex. This is immediate if $g \in \mathcal{N}$, since $\operatorname{val}(\mathcal{N}) \cap (\mathcal{G}/\mathcal{N})^{>} = \varnothing$. Suppose that $g \notin \mathcal{N}$. Then

$$
\begin{aligned}
\operatorname{val}(g\mathcal{N}) \cap (\mathcal{G}/\mathcal{N})^{>} &= \{f\mathcal{N} : f \asymp g\} \cap (\mathcal{G}/\mathcal{N})^{>} \\
&= \{f\mathcal{N} : f \asymp g \wedge f > 1\} \\
&= (\operatorname{val}(g) \cap \mathcal{G}^{>})\mathcal{N}.
\end{aligned}
$$

Since $\operatorname{val}(g) \cap \mathcal{G}^{>}$ is convex in $\mathcal{G}$, its image by the quotient map is convex in $\mathcal{G}/\mathcal{N}$, hence the result. Suppose now that $\mathcal{G}$ satisfies the growth axiom, and let $f, g \in \mathcal{G}$ with $f, g > \mathcal{N}$ and $f\mathcal{N} \succ g\mathcal{N}$. We have $f \succ g$ and $f, g > 1$ so $f\, g\, f^{-1} > g$ by **GA** in $\mathcal{G}$. Equivalently, $[f, g] < 1$, and since $f \not\asymp g$, we have $[f, g] \notin \mathcal{N}$, whence $[f, g] < \mathcal{N}$ by convexity of $\mathcal{N}$. So $[f, g]\mathcal{N} < \mathcal{N}$, i.e. **GA** holds in $\mathcal{G}/\mathcal{N}$. □

Given a $\rho \in \operatorname{val}(\mathcal{G})$, the group $\mathcal{C}_\rho = \mathcal{G}_{\leqslant \rho}/\mathcal{G}_{<\rho}$ of Proposition 2.19 is an ordered group for the natural ordering. For $f \in \mathcal{G}$ with $\rho = \operatorname{val}(f)$, we have $\operatorname{res}(f) = f\mathcal{G}_{<\rho} \in \mathcal{C}_\rho$. Therefore:

**Corollary 4.14.** *For all $f \in \mathcal{G}$, the set $\operatorname{res}(f)$ is convex.*

It follows that the ordering on each $\mathcal{C}_\rho$ extends to a total ordering $\lessdot$ on the set of residues $\{\operatorname{res}(g) : g \in \mathcal{G}\}$, given by

$$
\forall g, h \in \mathcal{G}, \operatorname{res}(g) \lessdot \operatorname{res}(h) \Longleftrightarrow \operatorname{res}(g) < \operatorname{res}(h).
$$

For $g, h \in \mathcal{G}$, we also write $g \lessdot h$ if $\operatorname{res}(g) \lessdot \operatorname{res}(h)$, so $\lessdot$ is a partial ordering on the set $\mathcal{G}$.

## 4.4 Growth order groups

We now introduce growth ordered groups. Let $(\mathcal{G}, \cdot, 1, <)$ be an ordered group satisfying the following $\mathcal{L}_{\mathrm{og}}$-sentence:

**GOG1.** For all $f, g \in \mathcal{G}^{>}$ with $f \geqslant g$ and all $g_0 \in \mathcal{C}(g)$, there is an $f_0 \in \mathcal{C}(f)$ with $f_0 \geqslant g_0$.

Define a binary relation $\preccurlyeq$ on $\mathcal{G}$ by setting, $f \preccurlyeq g$ for all $f, g \in \mathcal{G}$ such that $f$ lies in the convex hull of $\mathcal{C}(g)$. This is [8, Section 1.2] a linear quasi-ordering satisfying **D1**–**D5**. One sees in view of **D5** that it is the finest relation satisfying these axioms. Consider the following sentences in $\mathcal{L}_{\text{og}}$:

**GOG2.** For all $f, g \in \mathcal{G}^{>}$ with $f > \mathcal{C}(g)$, we have $f g > g f$.

**GOG3.** For all $f \in \mathcal{G}^{\neq}$, there is an $\mathfrak{s} \in \mathcal{G}^{>}$ such that $\mathcal{C}(\mathfrak{s})$ is Abelian, and for all $g \asymp f$, there is an $\mathfrak{s}_0 \in \mathcal{C}(\mathfrak{s})$ with $\mathfrak{s}_0 \sim \mathfrak{s}$.

So **GOG2** states that $(\mathcal{G}, \cdot, 1, <, \preccurlyeq)$ satisfies **GA** whereas **GOG3** states that $(\mathcal{G}, \cdot, 1, \preccurlyeq)$ has scaling elements. If $\mathcal{G}$ satisfies **GOG1** and **GOG3**, then [8, Lemma 2.18] the valued group $(\mathcal{G}, \cdot, 1, \preccurlyeq)$ satisfies **D6**. Therefore $(\mathcal{G}, \cdot, 1, <, \preccurlyeq)$ is an ordered c-valued group.

An ordered group satisfying **GOG1**, **GOG2** and **GOG3** is called a *growth order group*. In [8, Section 4], we showed that certain groups of germs at $+\infty$ of real-valued functions are growth order groups. In particular, if $\mathcal{R}$ is an o-minimal expansion (see [31, 14]) of the real ordered field which does not define the exponential function, then the group under composition of germs of germs at $+\infty$ of definable unary functions in $\mathcal{R}$ is a growth order group for the ordering of comparison at $+\infty$.

**Example 4.15.** Let $C$ be an ordered field and consider the ordered field $C((C))$ of generalised series with exponents and coefficients in $C$, i.e. of functions $s : C \longrightarrow C$ with anti-well-ordered support, with its natural valuation $v : C((C)) \longrightarrow C$ (see [6, Section 6]). This is an H-field for the standard derivation $\partial$ with $\partial(s)(c) = (c+1)\, s(c+1)$ for all $(s, c) \in C((C)) \times C$ Then the subset $\mathcal{P}_C \subset C((C))$ of parabolic series is an ordered group [6, Section 6.1]. Moreover, in view of [6, Section 6.2], it is a growth order group, whence an ordered c-valued group for the finest c-dominance relation. We also note that $\operatorname{val}(\mathcal{P}_C)$ is isomorphic to the underlying ordered set of $C$ and [6, Section 6.2] that each $\mathcal{C}_\rho$ for $\rho \in \operatorname{val}(\mathcal{P}_C)$ is isomorphic to $(C, +, 0, <)$.

**Example 4.16.** The set of positive infinite so-called finitely nested hyperseries of [5] is a growth order group [5, Section 9] for its natural composition law. It is an $\mathbb{R}$-valued $\mathbb{R}$-group for the $\mathbb{R}$-structure of [5, Section 9.2]. Any two positive elements in that ordered group are conjugate [5, Result E], so its underlying group is simple.

## 4.5 A coarser dominance relation

In this subsection, we fix an ordered valued group $(\mathcal{G}, \cdot, 1, <, \preccurlyeq)$ satisfying **GA**. We will define a coarser total quasi-ordering $\underline{\ll}$ on $\mathcal{G}$ which sometimes controls how far from being Abelian $\mathcal{G}$ is. For $f \in \mathcal{G}$, we write $f^{+} := \max(f, f^{-1})$ and $f^{-} := \min(f, f^{-1}) = (f^{+})^{-1}$.

**Definition 4.17.** *For $f, g \in \mathcal{G}^{\neq}$, we write*

$$g \,\underline{\ll}\, f \qquad \text{if} \qquad g f g^{-1} \sim f.$$

*We also write $g \ll f$ if $f \not\underline{\ll}\, g$, i.e. if $g f g^{-1} \nsim f$, or equivalently $[g^{-1}, f^{-1}] \succcurlyeq f$.*

Proposition 2.4 gives $g \preccurlyeq f \Longrightarrow g \,\underline{\ll}\, f$. Note that we have

$$g \,\underline{\ll}\, f \Longleftrightarrow [g^{-1}, f^{-1}] \prec f.$$

**Lemma 4.18.** *For all $f, g \in \mathcal{G}^{\neq}$, we have $g \ll f$ if and only if $f^{+} g^{+} f^{-} \gtrdot g^{+}$.*

**Proof.** Assume that $g \ll f$. So $f^+ \succ g^+$. Since $f^+ g^+ f^- > g^+$ by **GA** and $f^+ g f^- \nsim g$, we must have $f^+ g^+ f^- \gtrdot g^+$. Conversely, assume that $f^+ g^+ f^- \gtrdot g^+$, so $f^+ g^+ f^- \nsim g^+$. By taking inverses and using **D3**, we obtain $f^+ g f^- \nsim g$. Again with **D3**, by conjugating by $f^-$, we obtain $f g f^{-1} \nsim g$. □

**Lemma 4.19.** *Let $f, g, h \in \mathcal{G}$ with $f \preccurlyeq h$. If $g \ll f$, then $g \ll h$.*

**Proof.** We may assume that $g > 1$. Suppose first that $h^+ \nsim f^+$, so $h^+ f^- \asymp h$. Lemma 2.5 gives $h^+ f^- \succ f^+ g f^-$, so **GA** yields

$$h^+ g h^- = (h^+ f^-)(f^+ g f^-)(h^+ f^-)^{-1} > f^+ g f^-.$$

But $f^+ g f^- \gtrdot g$ by Lemma 4.18, so $h^+ g h^- \gtrdot g$, i.e. $g \ll h$ by Lemma 4.18.

Suppose now that $h^+ \sim f^+$ and assume for contradiction that $h g h^{-1} \sim g$. Conjugating by $h$, we obtain $h^2 g h^{-2} \sim h g h^{-1} \sim g$. But $(h^2)^+ \nsim f^+$, so our previous arguments and Lemma 4.18 entail that $h^2 g h^{-2} \nsim g$: a contradiction. So $g \ll h$. □

**Lemma 4.20.** *Let $f, g, h \in \mathcal{G}$ with $h \preccurlyeq g$. If $g \ll f$, then $h \ll f$. If $g \underline{\ll} f$, then $h \underline{\ll} f$.*

**Proof.** Assume for contradiction that $f h f^{-1} \sim h$ and set $g_0 := (f h^{-1}) g (f h^{-1})^{-1}$. We have

$$\begin{aligned} f g f^{-1} &= (f h^{-1})(h g h^{-1})(f h^{-1})^{-1} \\ &= (f h f^{-1}) g_0 (f h f^{-1})^{-1}. \end{aligned}$$

Now $f h^{-1} \asymp f$ so Lemma 4.19 gives $g \ll f h^{-1}$, so $g_0 \nsim g$. We have $f h f^{-1} \asymp h \nsucc g \asymp g_0$ so $g_0 \not\ll f h f^{-1}$, so $(f h f^{-1}) g_0 (f h f^{-1})^{-1} \sim g_0$. We deduce that $g \sim g_0$: a contradiction. □

**Lemma 4.21.** *Let $f, g, h \in \mathcal{G}$ with $h \preccurlyeq g$. If $g \underline{\ll} f$, then $h \underline{\ll} f$.*

**Proof.** Suppose that $g \underline{\ll} f$. We cannot have $f \ll h$, otherwise Lemma 4.19 would give $f \ll g$. So $h \underline{\ll} f$. □

**Lemma 4.22.** *Let $f, g, h \in \mathcal{G}$ with $f \preccurlyeq h$. If $g \underline{\ll} f$, then $g \underline{\ll} h$.*

**Proof.** Suppose that $g \underline{\ll} f$. We cannot have $h \ll g$, otherwise Lemma 4.20 would give $f \ll g$. So $g \underline{\ll} h$. □

Since $g \ll f \Longrightarrow g \preccurlyeq f$ for all $f, g \in \mathcal{G}$, we deduce that the relation $\ll$ is transitive. Thus:

**Corollary 4.23.** *The relation $\ll$ is a strict ordering on $\mathcal{G}^{\neq}$.*

**Definition 4.24.** *For all $f, g \in \mathcal{G}$, we write $g \prec\!\!\!\ll f$ if there is an $h \in \mathcal{G}^{\neq}$ with $g \ll h$ and $h \underline{\ll} f$. We write $g \underline{\prec\!\!\!\ll} f$ if $f \not\prec\!\!\!\ll g$. Lastly, we write $f \mathrel{\breve{\asymp}} g$ if $f \not\prec\!\!\!\ll g \wedge g \not\prec\!\!\!\ll f$.*

As direct corollaries of Lemmas 4.19, 4.20, 4.21 and 4.22, we have:

**Lemma 4.25.** *For all $f, g, h \in \mathcal{G}$ with $f \preccurlyeq h$, we have $g \prec\!\!\!\ll f \Longrightarrow g \prec\!\!\!\ll h$ and $g \underline{\prec\!\!\!\ll} f \Longrightarrow g \underline{\prec\!\!\!\ll} h$.*

**Lemma 4.26.** *For all $f, g, h \in \mathcal{G}$ with $g \prec\!\!\!\ll f$ and $h \preccurlyeq g$, we have $h \prec\!\!\!\ll f$.*

**Lemma 4.27.** *Let $f, g \in \mathcal{G}$. If $g \prec\!\!\!\ll f$, then $g \prec f$.*

**Proof.** Assume for contradiction that $g \prec\!\!\!\ll f$ and $f \preccurlyeq g$. Then $g \prec\!\!\!\ll g$ by Lemma 4.25, so there is an $h \in \mathcal{G}$ with $g \ll g$ and $h \underline{\ll} h$: a contradiction. □

**Proposition 4.28.** *The relation $\underline{\prec\!\!\!\ll}$ is a total quasi-ordering on $\mathcal{G}$ which is coarser than $\preccurlyeq$.*

**Proof.** Let $f, g, h \in \mathcal{G}$ with $f \underline{\prec\!\!\prec} g$ and $g \underline{\prec\!\!\prec} h$. Assume for contradiction that $f \not\underline{\prec\!\!\prec} h$, so $h \prec\!\!\prec f$. Let $j \in \mathcal{G}$ with $h \ll j$ and $j \underline{\ll} f$. We cannot have $j \underline{\ll} g$ since $h \not\prec\!\!\prec g$. So $g \ll j$, whence $g \prec\!\!\prec f$: a contradiction. This shows that $\underline{\prec\!\!\prec}$ is transitive. It is reflexive by Lemma 4.27 (and because $\prec$ is anti-reflexive). Let $f, g \in \mathcal{G}$, with $f \preccurlyeq g$. We cannot have $g \prec\!\!\prec f$ by Lemma 4.27, so $f \underline{\prec\!\!\prec} g$. Therefore $\underline{\prec\!\!\prec}$ is coarser than $\preccurlyeq$. We deduce since $\preccurlyeq$ is total that $\underline{\prec\!\!\prec}$ is total. □

**Remark 4.29.** We deduce (see Remark 2.2) that the total quasi-ordering $\underline{\prec\!\!\prec}$ on $\mathcal{G}$ satisfies both **D1**, **D2** and **D4**, and that it satisfies **D5** if $\preccurlyeq$ does. It also satisfies **D3** since conjugacy maps are automorphisms of ordered valued groups and $\underline{\prec\!\!\prec}$ is definable in $(\mathcal{G}, \cdot, 1, <, \preccurlyeq)$. The growth axiom holds, again because $\underline{\prec\!\!\prec}$ is coarser than $\preccurlyeq$. As for **D6**, it fails in general.

Let us illustrate what the relation $\prec\!\!\prec$ amounts to for the ordered groups of Section 4.1. Let $(K, +, \cdot, 0, 1, \mathcal{O}, \partial, \circ, x)$ be an H-field with composition, and denote by $\mathcal{o}$ the maximal ideal of $\mathcal{O}$ and write $v$ for the valuation corresponding to $\mathcal{O}$, i.e. the quotient map $K^{\times} \twoheadrightarrow K^{\times} / \mathcal{O}^{\times}$. Suppose that the ordered valued group $\mathcal{P} = \{x + \delta : \delta \in x\,\mathcal{o}\}$ of Proposition 4.8 satisfies **GA**. This is the case for instance if $K$ is the field of germs of definable maps in a levelled o-minimal expansion of the real ordered field [8, Theorem 2], or if $K$ is a substructure of the field of finitely nested hyperseries of Example 4.16.

**Lemma 4.30.** *Let all $\delta, \varepsilon \in x\,\mathcal{o} \setminus \{0\}$. If $\varepsilon \notin \mathcal{O}^{\times}$, then we have*

$$x + \delta \gg x + \varepsilon \Longleftrightarrow \varepsilon^{\dagger} \delta \notin \mathcal{o}.$$

*If $\varepsilon \in \mathcal{O}^{\times}$, then we have $x + \delta \underline{\ll} x + \varepsilon$.*

**Proof.** If $\varepsilon \notin \mathcal{O}^{\times}$, then $\max(\varepsilon, -\varepsilon) \in K^{>C}$ or $\max(\varepsilon^{-1}, -\varepsilon^{-1}) \in K^{>C}$. Thus this follows from Lemmas 4.6 and 4.7, and (4.2). Suppose now that $\varepsilon \in \mathcal{O}^{\times}$. Since $K$ is an H-field, we have $\varepsilon = c + \iota$ for a $\iota \in \mathcal{o}$ and a $c \in C^{\times}$. For each $\mu \in x\,\mathcal{o}$, we have $\varepsilon \circ (x + \mu) - \varepsilon = \iota \circ (x + \mu) - \iota$, so $\varepsilon \circ (x + \mu) - \varepsilon \in \varepsilon\,\mathcal{o} = \mathcal{o}$ if and only if $\iota \circ (x + \mu) - \iota \in \mathcal{o}$. But $\iota \in \mathcal{o}$ so $\iota \circ (x + \mu), \iota \in \mathcal{o}$ by [8, Section 4.1, **HFC1**]. So $\varepsilon \circ (x + \mu) - \varepsilon \in \varepsilon\,\mathcal{o}$ for all $\mu \in x\,\mathcal{o}$. We deduce by Lemma 4.7 and (4.2) that $(x + \delta) \circ (x + \varepsilon) \circ (x + \delta)^{\mathrm{inv}} - (x + \varepsilon) \in \varepsilon\,\mathcal{o}$, so $x + \delta \underline{\ll} x + \varepsilon$. □

**Corollary 4.31.** *For all $\delta, \varepsilon \in x\,\mathcal{o} \setminus \{0\}$, we have $x + \delta \not\underline{\succ\!\!\succ} x + \varepsilon$ if and only if*

- *$\delta \in \mathcal{O}^{\times}$, or*
- *$\delta, \varepsilon \in \mathcal{o}$ and there is no $\iota \in K$ with $v(\iota^{\dagger}) \in (v(\varepsilon), -v(\delta)]$.*

**Proof.** Suppose that $\delta \in \mathcal{O}^{\times}$. Then by Lemma 4.30 there is no $\iota \in x\,\mathcal{o}$ with $x + \delta \ll x + \iota$, whence in particular $x + \delta \not\underline{\succ\!\!\succ} x + \varepsilon$. We deduce by Lemma 4.25 that $x + \delta \not\underline{\succ\!\!\succ} x + \varepsilon$ whenever $\delta \notin \mathcal{o}$. Suppose that $\delta \in \mathcal{o}$. We have $x + \delta \not\underline{\succ\!\!\succ} x + \varepsilon$ if and only if there is no $\iota \in x\,\mathcal{o}$ with $x + \delta \ll x + \iota$ and $x + \iota \underline{\ll} x + \varepsilon$. By Lemma 4.30, this is equivalent to the conjunction that there be no $\iota \in x\,\mathcal{o}$ with $\iota^{\dagger} \delta \notin \mathcal{o} \wedge \varepsilon^{\dagger} \iota \in \mathcal{o}$. This translates into the statement above. □

In particular, if $K$ is Liouville-closed (see [2, p 3]), which entails that $v((K^{\times})^{\dagger}) = v(K)$, then the second case amounts to the conjunction that $\delta, \varepsilon \in \mathcal{o}$ and $v(\delta^{\dagger}) \geqslant v(\varepsilon^{\dagger})$. We will come bck to these matters in Section 6.

# 5 Valuative uniform structure and topology

We now consider the uniform structure on valued groups and related properties: valuative balls, Cauchy and pseudo-Cauchy sequences, completeness and spherical completeness. Our definitions are the same as in [35, Sections 2, 4 and 5] and [34], but we focus on the interaction of the uniform structure with the other axioms of c-valuations and $A$-valuations.

## 5.1 Topology and uniform structure

We fix a non-trivial valued group $(\mathcal{G},\cdot,1,\preccurlyeq)$. Given $g\in\mathcal{G}$ and $\rho\in\mathrm{val}(\mathcal{G})$, the (strict) *right ball* centered in $g$ with radius $\rho$ is the set

$$B_R(g,\rho):=\{h\in\mathcal{G}:\mathrm{val}(g^{-1}h)<\rho\}.$$

The (strict) *left ball* centered in $g$ with radius $\rho$ is the subset

$$B_L(g,\rho):=\{h\in\mathcal{G}:\mathrm{val}(h\,g^{-1})<\rho\}.$$

Note that $g\in B_L(g,\rho)\cap B_R(g,\rho)$. By **D3**, we have

$$B_R(g,\rho)=B_L(g,g\,\rho\,g^{-1}). \tag{5.1}$$

**Lemma 5.1.** *For all $g\in\mathcal{G}$, $\rho\in\mathrm{val}(\mathcal{G})$ and $h\in B_R(g,\rho)$, we have $B_R(h,\rho)=B_R(g,\rho)$.*

**Proof.** Write $\rho=\mathrm{val}(g_0)$ and let $f\in\mathcal{G}$. We have

$$h^{-1}f=(h^{-1}g)\,(g^{-1}f)$$

where $h^{-1}g\prec g_0$. By Proposition 2.4, we deduce that $h^{-1}f\prec g_0\Longleftrightarrow g^{-1}f\prec g_0$, hence the result. □

**Corollary 5.2.** *Let $f_1,f_2\in\mathcal{G}$ and $\rho_1,\rho_2\in\mathrm{val}(\mathcal{G})$. Then*

$$B_R(f_1,\rho_1)\subseteq B_R(f_2,\rho_2)\qquad\text{or}\qquad B_R(f_1,\rho_1)\cap B_R(f_2,\rho_2)=\varnothing.$$

By (5.1) and Corollary 5.2, the sets of right or left balls form a basis of neighborhoods of 1 in $\mathcal{G}$. By **D3**, for all $\delta\in\mathcal{G}^{\neq}$ and $f\in\mathcal{G}$, there is an an $\varepsilon\in\mathcal{G}^{\neq}$ such that for all $g\in\mathcal{G}$, $g\preccurlyeq\varepsilon\Longrightarrow f\,g\,f^{-1}\preccurlyeq\delta$ (take $\varepsilon=f^{-1}\,\delta\,f$). Therefore [35, Theorem 2.4] the resulting Hausdorff topological space is a topological group. We call that topology on $\mathcal{G}$ the *valuative topology*.

We have two uniform structures (as per [12, Section II.1]) on $\mathcal{G}$ whose entourages are supersets of sets $\{(f,g)\in\mathcal{G}\times\mathcal{G}:f\,g^{-1}\preccurlyeq h\}$ (left uniform structure) and $\{(f,g)\in\mathcal{G}\times\mathcal{G}:g^{-1}f\preccurlyeq h\}$ (right uniform structure) respectively, where $h$ ranges in $\mathrm{val}(\mathcal{G})$. By **D3**, those uniform structures coincide (see [35, Theorems 2.4 and 2.8]). We call that uniform structure the *valuative uniform structure* on $(\mathcal{G},\cdot,1,\preccurlyeq)$.

The product and inverse operations on $\mathcal{G}$ are uniformly continuous. Indeed this is clear for the inverse operation by **D4**. Given $(f_0,g_0,f_1,g_1)\in\mathcal{G}^4$ and $h_0,h_1,h_2\in\mathcal{G}^{\neq}$ with $f_0\,g_0^{-1}\preccurlyeq h_0$ and $f_1\,g_1^{-1}\preccurlyeq h_1$, we have

$$f_0\,f_1\,(g_0\,g_1)^{-1}=f_0\,(f_1\,g_1^{-1})\,f_0^{-1}\,f_0\,g_0^{-1}\preccurlyeq h_2$$

whenever $h_1\preccurlyeq f_0^{-1}\,h_2\,f_0$ and $h_0\preccurlyeq h_2$ by **D3**. So the product is uniformly continuous.

**Proposition 5.3.** *Let $A$ be a ring and let $(\mathcal{G},\cdot,1,(g\mapsto g^a)_{a\in A},\preccurlyeq)$ be a uniformly $A$-valued $A$-group. Then each $a$-power map $\mathcal{G}\longrightarrow\mathcal{G}\,;\,g\mapsto g^a$ for $a\in A$ is uniformly continuous.*

**Proof.** This follows immediately from $\mathbf{D9}^A$. □

## 5.2 Completeness and Cauchy sequences

We say that an embedding $\Psi:\mathcal{H}\longrightarrow\mathcal{G}$ between such valued groups is *coinitial* if $\mathrm{val}(\Phi(\mathcal{H}))$ is coinitial in $\mathrm{val}(\mathcal{G})$. In that case, the function $\Psi$ is uniformly continuous. Indeed, for all $f,g\in\mathcal{H}$ and $\varepsilon\in\mathcal{G}^{\neq}$, picking a $\delta\in\mathcal{H}$ with $\Phi(\delta)\preccurlyeq\varepsilon$, we have $\Psi(f)\,\Psi(g)^{-1}\preccurlyeq\varepsilon$ whenever $f\,g^{-1}\preccurlyeq\delta$.

An extension $\Psi:\mathcal{H}\longrightarrow\mathcal{G}$ of valued groups is said *dense* if $\Psi(\mathcal{H})$ is dense in $\mathcal{G}$. Note that if the extension is dense, then $\mathrm{val}(\Psi(\mathcal{H}))=\mathrm{val}(\mathcal{G})$. In particular, dense extensions are initial.

**Lemma 5.4.** *Let $\Psi:\mathcal{H}\longrightarrow\mathcal{G}$ be a dense extension of valued groups. Then for each $\rho\in\mathrm{val}(\mathcal{G})$, the embedding $\mathrm{res}_\Psi:\mathcal{C}_{\mathrm{val}_\Phi^{-1}(\rho)}\longrightarrow\mathcal{C}_\rho$ is an isomorphism.*

**Proof.** Let $g\in\mathcal{G}$ with $\mathrm{val}(g)=\mathrm{val}_\Phi(\rho)$. By density, there is an $h\in\mathcal{H}$ with $g\sim\Phi(h)$ so $\mathrm{res}(\Psi(h))=\mathrm{res}(g)$. □

**Proposition 5.5.** *Let $A$ be a ring. All dense extensions of an $A$-torsion-free uniformly $A$-valued $A$-group are $A$-torsion-free.*

**Proof.** This follows from Lemmas 3.12 and 5.4. □

**Proposition 5.6.** *Let $\Phi:\mathcal{H}\longrightarrow\mathcal{G}$ be a dense extension of valued groups. If $\flat\in\mathrm{val}(\mathcal{H})$ is a commutator bound for an $(f,g)\in\mathcal{H}^2$, then $\mathrm{val}_\Phi(\flat)$ is a commutator bound for $(\Phi(f),\Phi(g))$.*

**Proof.** Let $h,j\in\mathcal{G}$ with $h\sim\Phi(f)$ and $j\sim\Phi(g)$. By density of the extension, there is an $(f_0,g_0)\in\mathcal{H}^2$ such that $\varepsilon:=\Phi(f_0)^{-1}h$ and $\delta:=\Phi(g_0)^{-1}j$ satisfy $\mathrm{val}(\delta),\mathrm{val}(\varepsilon)<\mathrm{val}_\Phi(\flat)$. In particular $f_0\sim f$ and $g_0\sim g$, so $[f_0,g_0]\succcurlyeq\flat$. Now

$$\begin{aligned}
[h,j] &= [\Phi(f_0)\,\varepsilon,\Phi(g_0)\,\delta]\\
&= \varepsilon^{-1}\,\Phi(f_0)^{-1}\,\delta^{-1}\,\Phi(g_0)^{-1}\,\Phi(f_0)\,\varepsilon\,\Phi(g_0)\,\delta\\
&= \varepsilon^{-1}\,(\Phi(f_0)^{-1}\,\delta^{-1}\,\Phi(f_0))\,\Phi(f_0)^{-1}\,\Phi(g_0)^{-1}\,\Phi(f_0)\,\Phi(g_0)\,(\Phi(g_0)^{-1}\,\varepsilon\,\Phi(g_0))\,\delta\\
&= \varepsilon^{-1}\,(\Phi(f_0)^{-1}\,\delta^{-1}\,\Phi(f_0))\,[\Phi(f_0),\Phi(g_0)]\,(\Phi(g_0)^{-1}\,\varepsilon\,\Phi(g_0))\,\delta.
\end{aligned}$$

Now $\mathrm{val}([\Phi(f_0),\Phi(g_0)])\geqslant\mathrm{val}_\Phi(\flat)$. By continuity of the maps $x\mapsto\Phi(f_0)^{-1}\,x^{-1}\,\Phi(f_0)$ and $y\mapsto\Phi(g_0)^{-1}\,y\,\Phi(g_0)$ and uniformity of the commutator bound $\flat$ in $\mathrm{res}(f_0),\mathrm{res}(g_0)$ we may choose $f_0,g_0$ such that $\varepsilon,\delta$ are small enough in valuation, so that $\mathrm{val}(\Phi(f_0)^{-1}\,\delta^{-1}\,\Phi(f_0))$, $\mathrm{val}(\Phi(g_0)^{-1}\,\varepsilon\,\Phi(g_0))<\mathrm{val}_\Phi(\flat)$. We then obtain $[h,j]\asymp[\Phi(f_0),\Phi(g_0)]$ by Proposition 2.4, so $\mathrm{val}([h,j])\geqslant\mathrm{val}_\Phi(\flat)$. □

We wee with Lemma 2.29 that:

**Corollary 5.7.** *Let $\Phi:\mathcal{H}\longrightarrow\mathcal{G}$ be a dense extension of valued groups where $\mathcal{H}$ is c-bounded. Then $\mathcal{G}$ is c-bounded. In particular, it satisfies* **D5**.

**Proposition 5.8.** *Let $\Phi:\mathcal{H}\longrightarrow\mathcal{G}$ be a dense extension of valued groups where $\mathcal{H}$ is c-bounded. If $\mathcal{H}$ is c-valued, then so is $\mathcal{G}$.*

**Proof.** It suffices to show that **D6** holds in $\mathcal{G}$. This is immediate, by continuity of the commutator map. □

**Definition 5.9.** *Let $(\mathcal{G},\cdot,1,\preccurlyeq)$ be a valued group. A* **Cauchy sequence** *in $\mathcal{G}$ is a sequence $(g_\gamma)_{\gamma<\lambda}\in\mathcal{G}^\lambda$ where $\lambda$ is a non-zero limit ordinal such that for all $\varepsilon\in\mathcal{G}^{\neq}$, there is a $\mu<\lambda$ such that for all $\gamma,\eta\in[\mu,\lambda)$, we have $g_\gamma\,g_\eta^{-1}\prec\varepsilon$.*

A sequence $(g_\gamma)_{\gamma<\lambda}\in\mathcal{G}^\lambda$ for non-zero limit ordinal $\lambda$ converges to a $g\in\mathcal{G}$ if it converges to $g$ at $+\infty$ as a function $\lambda\longrightarrow\mathcal{G}$ for the order topology on $\lambda$, i.e. if for all $\varepsilon\in\mathcal{G}^{\neq}$, there is a $\mu<\lambda$ such that $g_\gamma g^{-1}\prec\varepsilon$ whenever $\gamma\in[\mu,\lambda)$. By **D3**-**D4**, one can replace that condition by $\forall\varepsilon\in\mathcal{G}^{\neq},\exists\mu<\lambda,\forall\gamma\in[\mu,\lambda),g^{-1}\,g_\gamma\prec\varepsilon$. It is easy to see that convergent sequences are Cauchy sequences (see [35, Proposition 4.2]).

**Definition 5.10.** *A valued group $\mathcal{G}$ is said* **complete** *if its Cauchy sequences are convergent.*

**Theorem 5.11.** *Let $\mathcal{G}$ be a valued group. There is a complete valued group $\tilde{\mathcal{G}}$ and a dense embedding $\Psi:\mathcal{G}\longrightarrow\tilde{\mathcal{G}}$ such that for any dense extension $\Phi:\mathcal{G}\longrightarrow\mathcal{H}$, there is a unique $\Lambda:\mathcal{H}\longrightarrow\tilde{\mathcal{G}}$ with $\Psi=\Lambda\circ\Phi$. Moreover $\tilde{\mathcal{G}}$ is c-bounded if $\mathcal{G}$ is c-bounded, and it is c-valued if $\mathcal{G}$ is c-bounded and c-valued.*

**Proof.** Consider the completion $\Psi:\mathcal{G}\longrightarrow\tilde{\mathcal{G}}$ of $\mathcal{G}$ as a uniform group [12, Section III.4, Theorem 1]. This is a dense extension of uniform groups which is final among dense extensions of uniform groups of $\mathcal{G}$ [12, Section II.7, Theorem 3 and Proposition 12]. Therefore it suffices to endow $\tilde{\mathcal{G}}$ with a dominance relation such that $\Psi$ is an embedding of valued groups. We define $\tilde{\preccurlyeq}$ to be the closure in $\tilde{\mathcal{G}}$ of the set $\{(\Psi(f),\Psi(g)):(f,g)\in\mathcal{G}^2\wedge f\preccurlyeq g\}$. By continuity of the product, conjugacy and inverse maps on $\mathcal{G}\times\mathcal{G}$ and $\mathcal{G}$ respectively, the axioms **D1–D4** hold in $(\tilde{\mathcal{G}},\cdot,1,\tilde{\preccurlyeq})$, and $\Psi$ is an embedding by definition. The remaining parts of the statement of the theorem follow from Corollary 5.7 and Proposition 5.8. □

We call $\tilde{\mathcal{G}}$ the *completion* of $\mathcal{G}$. It can also be obtained directly [34, Theorem 4.4] as the quotient of the set of Cauchy sequences in $\mathcal{G}$ indexed by the cofinality of $(\mathrm{val}(\mathcal{G}),>)$ by its normal subgroup of Cauchy sequences that converge to 1. In the case when each subgroup $\mathcal{G}_{<\rho}$, $\rho\in\mathrm{val}(\mathcal{G})$ is normal, it can be obtained as the inverse limit of the inverse system of groups $(\mathcal{G}/\mathcal{G}_{<\rho})_{\rho\in\mathrm{val}(\mathcal{G})}$. For instance, completions of non-Abelian free exponential groups over commutative domains of characteristic 0 coincide with their pronilpotent completions.

Recall by [12, Section II.6, Theorem 2] that any uniformly continuous function on a valued group $\mathcal{G}$ extends uniquely into a uniformly continuous function on $\tilde{\mathcal{G}}$.

**Proposition 5.12.** *Let $A$ be a ring and let $\mathcal{G}$ be a uniformly $A$-valued $A$-group. Then $\tilde{\mathcal{G}}$ is a uniformly $A$-valued $A$-group for the uniformly continuous extension of all $a$-power maps $\mathcal{G}\longrightarrow\mathcal{G}$.*

**Proof.** This follows from the above comment and Proposition 5.3, since the axioms for $A$-groups and $A$-valued $A$-groups are topologically closed. □

**Corollary 5.13.** *Let $A$ be a unital $\mathbb{Q}$-algebra and let $\mathcal{G}$ be a uniformly $A$-valued $A$-group. Then $\tilde{\mathcal{G}}$ is divisible.*

**Proposition 5.14.** *Let $A$ be a ring, let $\beta\leqslant\omega$ and let $(\mathcal{G}_{[n]})_{n\leqslant\beta}$ be an $A$-series with associated dominance relation $\preccurlyeq$ on $\mathcal{G}=\mathcal{G}_{[0]}$. Let $(\tilde{\mathcal{G}},\cdot,1,\tilde{\preccurlyeq})$ be the completion of $\mathcal{G}$. For each $n\leqslant\omega$, write $\tilde{\mathcal{G}}_{[n]}$ for the relative closure of $\mathcal{G}_{[n]}$ in $\tilde{\mathcal{G}}$. Then $(\tilde{\mathcal{G}}_{[n]})_{n\leqslant\beta}$ is an $A$-series with associated dominance relation $\tilde{\preccurlyeq}$.*

**Proof.** If $\beta<\omega$, then $\tilde{\mathcal{G}}_{[\beta]}$ is the topological closure of $\{1\}$, so $\tilde{\mathcal{G}}_{[\beta]}=\{1\}$. If $\beta=\omega$ then that $\bigcap_{n<\beta}\tilde{\mathcal{G}}_{[n]}\subseteq\{1\}$ follows from the fact that $\widetilde{\mathrm{val}}(\bigcap_{m\leqslant n}\tilde{\mathcal{G}}_{[m]})\geqslant n$ in $(\omega,>)$ for all $n<\omega$. The other axioms of $A$-series are seen to hold by uniform continuity of the extended $a$-power maps and the group operations. We have $\tilde{\mathcal{G}}_{[0]}=\tilde{\mathcal{G}}$ by density of $\mathcal{G}$ in $\tilde{\mathcal{G}}$. Lastly, by definition of the dominance relation $\tilde{\preccurlyeq}$ on $\tilde{\mathcal{G}}$, we see that it coincides with that associated to the $A$-series. □

**Remark 5.15.** This mail fail for $A$-series of length $>\omega$, for the topological closure of a decreasing intersection may be strictly smaller than the intersection of topological closures.

**Remark 5.16.** We note by uniform continuity of commutators that $\tilde{\mathcal{G}}_{[n]} = \tilde{\mathcal{G}}_n$ (resp. $\tilde{\mathcal{G}}_{[n]} = \tilde{\mathcal{G}}^{(n)}$) for all $n \leqslant \beta$ if $(\mathcal{G}_{[n]})_{n \leqslant \beta}$ is the central lower (resp. derived) $A$-series on $\mathcal{G}$.

**Corollary 5.17.** *Let $A$ be a torsion-free ring and let $\mathcal{H}$ be a residually $A$-nilpotent (resp. residually $A$-solvable) $A$-group with its lower central (resp. derived) valuation. Then $\tilde{\mathcal{H}}$ is a residually $A$-nilpotent (resp. residually $A$-solvable) $A$-group.*

In the case when $\mathcal{G}$ is a c-bounded ordered valued group, one extends the ordering on $\mathcal{G}$ to an ordering on the completion $\tilde{\mathcal{G}}$ of its underlying valued group by density, and obtains:

**Theorem 5.18.** *Let $\mathcal{G}$ be a c-bounded ordered valued group. There is a complete, c-bounded ordered valued group $\tilde{\mathcal{G}}$ and a dense embedding $\Psi: \mathcal{G} \longrightarrow \tilde{\mathcal{G}}$ such that for any dense extension $\Phi: \mathcal{G} \longrightarrow \mathcal{H}$, there is a unique $\Lambda: \mathcal{H} \longrightarrow \tilde{\mathcal{G}}$ with $\Psi = \Lambda \circ \Phi$. Moreover $\tilde{\mathcal{G}}$ is c-valued if $\mathcal{G}$ is c-valued.*

**Proposition 5.19.** *Let $\mathcal{G}$ be a divisible ordered c-valued group. All dense and complete extensions of $\mathcal{G}$ are divisible.*

**Proof.** Recall that $\mathcal{G}$ is torsion-free, whence by Lemma 2.14 each element of $\mathcal{G}$ has a unique $n$-th root for all $n > 0$. It suffices to show that each $n$-th root map for $n > 0$ is uniformly continuous on $\mathcal{G}$. Let $\delta \in \mathcal{G}^{\neq}$. Let $f, g \in \mathcal{G}$ and write $\varepsilon := f g^{-1}$. We have $f^n g^{-n} = \prod_{i \leqslant n} f^i \varepsilon f^{-i}$. Since each $f^i \varepsilon f^{-i}$ for $i \leqslant n$ has the same sign as $\varepsilon$, it follows that $f^n g^{-n} \succsim \varepsilon$. This implies that the $n$-th root map is uniformly continuous. □

## 5.3 Spherical completeness and pseudo-Cauchy sequences

Let $(\mathcal{G}, \cdot, 1, \preccurlyeq)$ be a valued group.

**Definition 5.20.** *We say that $\mathcal{G}$ is* **spherically complete** *if any decreasing family of right and left balls has non-empty intersection.*

In view of (5.1) and Corollary 5.2, we have:

**Proposition 5.21.** *The following assertions are equivalent:*

*a) $\mathcal{G}$ is spherically complete;*

*b) any decreasing family of right (resp. left) balls has non-empty intersection;*

*c) any family of balls with the finite intersection property has non-empty intersection;*

*d) any family of right (resp. left) balls with the finite intersection property has non-empty intersection.*

**Definition 5.22.** *Let $\lambda$ be a non-zero ordinal. A sequence $(g_\gamma)_{\gamma < \lambda} \in \mathcal{G}$ is said* **pseudo-Cauchy** *if for all $\gamma < \eta < \mu < \lambda$, we have*

$$g_\eta g_\gamma^{-1} \succ g_\mu g_\eta^{-1}. \tag{5.2}$$

Taking inverses and conjugating by $g_\eta^{-1}$, we see that it is equivalent to replace (5.2) by

$$g_\gamma^{-1} g_\eta \succ g_\eta^{-1} g_\mu. \tag{5.3}$$

**Definition 5.23.** *Let $g \in \mathcal{G}$. We say that a pseudo-Cauchy sequence $(g_\gamma)_{\gamma < \lambda}$* **pseudo-converges** *to $g$, or that $g$ is a* **pseudo-limit** *of $(g_\gamma)_{\gamma < \lambda}$ if $\lambda$ is a limit and the sequence $(v(g_\gamma g^{-1}))_{\gamma < \lambda}$ is strictly decreasing.*

Again, one may replace each $g_\gamma g^{-1}$ by $g^{-1} g_\gamma$ in the definition. It is easy to see that:

**Lemma 5.24.** *Given a non-zero limit ordinal $\lambda$, a pseudo-Cauchy sequence $(g_\gamma)_{\gamma<\lambda}$ in $\mathcal{G}$ and an $f\in\mathcal{G}$, the sequences $(f g_\gamma)_{\gamma<\lambda}$, $(g_\gamma f)_{\gamma<\lambda}$ and $(g_\gamma^{-1})_{\gamma<\lambda}$ are pseudo-Cauchy. If moreover $(g_\gamma)_{\gamma<\lambda}$ pseudo-converges to a $g\in\mathcal{G}$, then they pseudo-converge to $f g$, $g f$ and $g^{-1}$ respectively.*

**Lemma 5.25.** *Given a non-zero limit ordinal $\lambda$, a pseudo-Cauchy sequence $(g_\gamma)_{\gamma<\lambda}$ in $\mathcal{G}$ and an element $g\in\mathcal{G}$, we have $(g_\gamma)_{\gamma<\lambda}\rightsquigarrow g$ if and only if $\forall\gamma<\lambda, g g_\gamma^{-1}\asymp g_{\gamma+1} g_\gamma^{-1}$.*

**Proof.** The sequence $(\mathrm{val}(g_{\gamma+1} g_\gamma^{-1}))_{\gamma<\lambda}$ is strictly decreasing by (5.2), so it suffices to show the left to right implication. Suppose that $(g_\gamma)_{\gamma<\lambda}\rightsquigarrow g$ and let $\gamma<\lambda$. We have $g g_\gamma^{-1}\succ g g_{\gamma+1}^{-1}$. Now $g g_\gamma^{-1}(g g_{\gamma+1}^{-1})^{-1}= g_{\gamma+1} g_\gamma^{-1}$ so Proposition 2.4 yields $g_{\gamma+1} g_\gamma^{-1}\asymp g g_\gamma^{-1}$, as claimed. □

**Proposition 5.26.** *Let $\lambda$ be a non-zero limit ordinal and let $(g_\gamma)_{\gamma<\lambda}$ be a pseudo-Cauchy sequence in $\mathcal{G}$. Then $(g_\gamma)_{\gamma<\lambda}$ pseudo-converges in an elementary extension of $\mathcal{G}$*

**Proof.** Consider the elementary diagram $\Delta(\mathcal{G})$ of $\mathcal{G}$ in the expansion of $\mathcal{L}_{\mathrm{vg}}$ with constant symbols for elements in $\mathcal{G}$. Write $\mathcal{L}^*$ for the expansion of that language with one additional constant symbol $z$. Since $(g_\gamma)_{\gamma<\lambda}$ is pseudo-Cauchy, the theory $\Delta(\mathcal{G})\cup\{g_\eta g_\gamma^{-1}\succ z g_\eta^{-1}: \gamma<\eta<\lambda\}$ is finitely consistent: each finite subset is realised in $\mathcal{G}$ by interpreting $z$ by an element of the sequence of sufficiently large index. By compactness, this theory has a model $\mathcal{G}^*$. Let $g\in\mathcal{G}^*$ interpret $z$. We see with Lemma 5.25 that $g$ is a pseudo-limit of $(g_\gamma)_{\gamma<\lambda}$. □

**Proposition 5.27.** *The valued group $\mathcal{G}$ is spherically complete if and only if all its pseudo-Cauchy sequences have pseudo-limits in $\mathcal{G}$.*

**Proof.** The proof is the same as [4, Lemma 2.2.10]. □

# 6 Near Abelianness

## 6.1 Nearly Abelian valued groups

**Definition 6.1.** *A valued group $(\mathcal{G},\cdot,1,\preccurlyeq)$ is said* **nearly Abelian** *if for all $f,g\in\mathcal{G}^{\neq}$, we have*

$$[f,g]\prec f,g. \tag{6.1}$$

This is a strengthening of the axiom **D6**, as here no assumption is made that $f\asymp g$. Note that any Abelian valued group is nearly Abelian. It is immediate that:

**Proposition 6.2.** *If $\mathcal{G}$ is a nearly Abelian valued group and $\mathcal{H}\subseteq\mathcal{G}$ is a subgroup, then $\mathcal{H}$ is nearly Abelian for the induced quasi-ordering.*

**Remark 6.3.** If $(\mathcal{G},\cdot,1,\preccurlyeq)$ is a non-Abelian, nearly Abelian valued group, then $\mathrm{val}(\mathcal{G})$ has no minimum. Indeed, for $\rho=\mathrm{val}(f)\in\mathrm{val}(\mathcal{G})$, we find a $g\in\mathcal{G}$ with $g\asymp f$ and $[f,g]\neq 1$ since $\mathcal{G}$ is non-Abelian. We have $[f,g]\prec f$, so $\mathrm{val}([f,g])\in\mathrm{val}(\mathcal{G})$ lies strictly below $\rho$.

**Corollary 6.4.** *A nearly Abelian valued group contains no non-Abelian subgroup of finite value set.*

**Example 6.5.** As a consequence of Example 2.11(vi), the c-valued group $C((t))^{\equiv 1}$ of Example 2.11 is nearly Abelian.

**Proposition 6.6.** *Let $A$ be a ring and let $\mathcal{G}$ be an $A$-hypocentral $A$-group with its lower central valuation $v$. Then $\mathcal{G}$ is nearly Abelian.*

**Proof.** For $f, g \in \mathcal{G}^{\neq}$, we have $[f, g] \in \mathcal{G}_{\max(v(f), v(g))+1}$, hence the result. □

**Proposition 6.7.** *Let $(\mathcal{G}, \cdot, 1, \preccurlyeq)$ be a valued group. The following assertions are equivalent:*

*a) $\mathcal{G}$ is nearly Abelian.*

*b) For all $g \in \mathcal{G}^{\neq}$ and all $f \in \mathcal{G}$, the natural isomorphism $\mathcal{C}_{v(g)} \longrightarrow \mathcal{C}_{v(fgf^{-1})}$ is trivial.*

*c) For all $f, g \in \mathcal{G}$, we have $f g f^{-1} \sim g$.*

*d) For all $f, g \in \mathcal{G}$, we have $f g \sim g f$.*

**Proof.** (a ⇒ b) Suppose that a holds. For $f, g \in \mathcal{G}$, we have $f\,(g f^{-1}) \sim (g f^{-1})\, f \sim g$, whence b holds.

(b ⇔ c) It is clear that b and c are equivalent.

(c ⇒ d) Suppose that c holds. For $f, g \in \mathcal{G}$, we have $f g = f\,(g f)\, f^{-1} \sim g f$, so d holds.

(d ⇒ c) Suppose that d holds. For $f, g \in \mathcal{G}$, we have $g \sim (g f^{-1})\, f \sim f\,(g f^{-1})$, so c holds.

(c ⇒ a) Suppose that c holds and let $f, g \in \mathcal{G}^{\neq}$. We have $g^{-1} f^{-1} g \sim f^{-1}$ by c, so $[g, f] \prec f$, so $[f, g] = [g, f]^{-1} \prec f$. Therefore a holds. This concludes the proof. □

**Proposition 6.8.** *Any nearly Abelian valued $\mathbb{Z}$-group satisfying* **D7**$^A$ *is uniformly $\mathbb{Z}$-valued.*

**Proof.** Let $(\mathcal{G}, \cdot, 1, \preccurlyeq)$ be nearly Abelian and satisfy **D7**$^A$ for $A = \mathbb{Z}$. Let $f \in \mathcal{G}^{\neq}$ and $\varepsilon \in \mathcal{G}$ with $\varepsilon \prec f$, and let $n \in \mathbb{Z} \setminus \{0\}$. We have $(f\varepsilon)^{n+1} f^{-(n+1)} = f\,(\varepsilon\,(f\varepsilon)^n f^{-n})\, f^{-1} \sim \varepsilon\,(f\varepsilon)^n f^{-n}$ by Proposition 6.7. So $f \varepsilon f^{-1} \preccurlyeq \varepsilon$. We deduce by induction and by Proposition 2.4 that $(f\varepsilon)^n f^{-n} \preccurlyeq \varepsilon$, i.e. **D9**$^A$ holds. □

A subgroup $\mathcal{H}$ of a group $\mathcal{G}$ is said *malnormal* if $\mathcal{H} \cap g \mathcal{H} g^{-1} = \{1\}$ for all $g \in \mathcal{G} \setminus \mathcal{H}$. A *CSA-group* is a group in which all maximal Abelian subgroups are malnormal.

**Proposition 6.9.** *If $(\mathcal{G}, \cdot, 1, \preccurlyeq)$ is a nearly Abelian c-valued group, then $(\mathcal{G}, \cdot, 1)$ is a CSA-group.*

**Proof.** Let $f \in \mathcal{G}^{\neq}$. By Proposition 2.12, it suffices to show that $\mathcal{C}(f)$ is malnormal. For $g \in \mathcal{G} \setminus \mathcal{C}(f)$ and $h \in \mathcal{C}(f)^{\neq}$, we have $g h g^{-1} \sim h$ by Proposition 6.7. We also have $g h g^{-1} \neq h$. Since there cannot be two distinct elements of $\mathcal{C}(f)$ with the same residue, it follows that $g h g^{-1} \notin \mathcal{C}(f)$. □

**Lemma 6.10.** *Let $\Psi : \mathcal{H} \longrightarrow \mathcal{G}$ be a dense extension of valued groups. If $\mathcal{H}$ is nearly Abelian, then so is $\mathcal{G}$.*

**Proof.** For simplicity, we assume that $\mathcal{H} \subseteq \mathcal{G}$ and that $\Psi$ is the inclusion. If $\mathcal{H}$ is Abelian, then $\mathcal{H} = \mathcal{G}$. So we may assume that $\mathcal{H}$ is non-Abelian, whence by Remark 6.3, $\mathrm{val}(\mathcal{H})$ has no minimum. Let $f, g \in \mathcal{G}^{\neq}$. The commutator function $(x, y) \mapsto [x, y]$ is continuous , so for sufficiently small $\rho \in \mathrm{val}(\mathcal{G})$ with $\rho < \mathrm{val}(f), \mathrm{val}(g), \mathrm{val}(f^{-1} g)$, we have $[f', g'] \sim [f, g]^{-1}$ whenever $(f', g') \in B_L(f, \rho) \times B_L(g, \rho)$. By density, we find such a pair $(f', g')$ in $\mathcal{H}^2$. Our choice of $\rho$ and the fact that $\mathcal{H}$ is nearly Abelian entail that $[f', g'] \prec f', g'$ where $f' \asymp f$ and $g' \asymp g$. Therefore $[f, g] \prec f, g$. So $\mathcal{G}$ is nearly Abelian. □

**Proposition 6.11.** *Let* $(\mathcal{G},\cdot,1,<,\preccurlyeq)$ *be a nearly Abelian ordered valued group satisfying* **GA**. *Then* $(\mathcal{G},\cdot,1,\preccurlyeq)$ *is c-bounded.*

**Proof.** Let $f,g\in\mathcal{G}^{\neq}$ with $f\succ g$. $(\rho,\mu)=(\mathrm{val}(f),\mathrm{val}(g))$. Note that $[f^{-1},g^{-1}]=(g\,f)\,[f,g]\,(g\,f)^{-1}\sim[f,g]\sim f\,[f,g]\,f^{-1}=[f^{-1},g]$ and that $[f,g^{-1}]=f^{-1}[g^{-1},f]\,f\sim[g^{-1},f]\sim[f^{-1},g]\sim[f,g]$. So we may assume that $f$ and $g$ are both positive.

Pick an $h\in\mathcal{G}^{>}$ with $h\asymp g$ and an $\varepsilon\in\mathcal{G}^{>}$ with $\varepsilon\prec g$. We claim that $\flat:=\mathrm{val}([h^{-1},\varepsilon^{-1}])$ is a c-bound for $(f,g)$. Indeed, we have

$$\begin{aligned}
[f^{-1},g^{-1}] &\succcurlyeq [f^{-1},\varepsilon^{-1}], && \text{(by Lemma 4.12)}\\
[f^{-1},\varepsilon^{-1}] &\succcurlyeq [h^{-1},\varepsilon^{-1}] && \text{(by Lemma 4.11).}
\end{aligned}$$

Since $\flat$ only depends on $\mathrm{val}(g)$, this concludes the proof. ☐

## 6.2 Groups of contracting strongly linear derivations

Let $C$ be a field of characteristic zero and $(\Gamma,+,0,<)$ be a non-trivial ordered *Abelian* group. Consider the field $C((\Gamma))$ of Hahn series [18] with coefficients in $C$ and exponents in $\Gamma$. This is the ring, under pointwise sum and Cauchy product, of functions $f:\Gamma\longrightarrow C$ whose support $\mathrm{supp}\,f=\{\gamma\in\Gamma:f(\gamma)\neq 0\}$ is a well-ordered subset of $\Gamma$. We have a canonical inclusion $C\longrightarrow C((\Gamma))$ whereby a $c\in C$ is identified with the function with support $\{0\}$ and value $c$ at 0. We have a dominance relation [4, Definition 3.3.1] $\preccurlyeq$ on the Abelian group $(C((\Gamma)),+,0)$ given by $f\preccurlyeq g\Longleftrightarrow(\min\mathrm{supp}\,f)\geqslant(\min\mathrm{supp}\,g)$ for all $f,g\neq 0$, and $0\preccurlyeq h$ for all $h$. It corresponds to the canonical valuation on the field $C((\Gamma))$. A linear map $\phi:C((\Gamma))\longrightarrow C((\Gamma))$ is said *contracting* if $\phi(f)\prec f$ for all $f\neq 0$.

There is a notion of infinite sum for certain families in $C((\Gamma))$ called summable families (see [20]). A $C$-linear map $C((\Gamma))\longrightarrow C((\Gamma))$ which commutes with infinite sums of summable families is said *strongly linear*. Let $\partial:C((\Gamma))\longrightarrow C((\Gamma))$ be a fixed strongly linear derivation with kernel $\mathrm{Ker}(\partial)=C$. Write

$$\mathrm{Cont}(\partial):=\{f\in C((\Gamma)):f\,\partial\text{ is contracting}\}.$$

We have a Lie bracket $[\![\cdot,\cdot]\!]:\mathrm{Cont}(\partial)\times\mathrm{Cont}(\partial)\longrightarrow\mathrm{Cont}(\partial);(f,g)\mapsto f\,\partial(g)-\partial(f)\,g$ on $\mathrm{Cont}(\partial)$. In [7], we showed that $\mathrm{Cont}(\partial)$ is a group for the operation

$$f*g:=f+g+\frac{1}{2}[\![f,g]\!]+\frac{1}{12}([\![f,[\![f,g]\!]]\!]-[\![g,[\![f,g]\!]]\!])+\cdots. \tag{6.2}$$

**Lemma 6.12.** [7, Lemma 1.1] *For* $f,g\in\mathrm{Cont}(\partial)\backslash\{0\}$, *we have* $[\![f,g]\!]\prec f,g$.

**Lemma 6.13.** *For* $f,g\in\mathrm{Cont}(\partial)\backslash\{0\}$, *we have* $[\![f,g]\!]=0$ *if and only if* $f\in C^{\times}g$.

**Proof.** It suffices to show that $f$ and $g$ are linearly dependent if $[\![f,g]\!]=0$. Suppose that $[\![f,g]\!]=0$. This means that $f$ and $g$ are solutions of the linear differential equation $\partial(y)-\frac{\partial(g)}{g}\,y=0$. But we see that a quotient of two non-zero solutions of that equation must lie in $\mathrm{Ker}(\partial)=C$. ☐

**Lemma 6.14.** *For* $f,g\in\mathrm{Cont}(\partial)\setminus\{0\}$, *we have* $f*g=g*f$ *if and only if* $f\in C^{\times}g$.

**Proof.** Assume for contradiction that $f*g=g*f$ and $[\![f,g]\!]\neq 0$. Note that $0=f*g-g*f=[\![f,g]\!]+\frac{1}{12}([\![f,[\![f,g]\!]]\!]-[\![g,[\![f,g]\!]]\!]-[\![g,[\![g,f]\!]]\!]+[\![f,[\![g,f]\!]]\!])+\cdots$. We see with Lemma 6.12 that $0=[\![f,g]\!]+h$ for an $h\prec[\![f,g]\!]$: a contradiction. Therefore $[\![f,g]\!]=0$. We conclude with Lemma 6.13. ☐

**Proposition 6.15.** *The structure* $(\mathrm{Cont}(\partial), *, 0, \preccurlyeq)$ *is a nearly Abelian c-valued group.*

**Proof.** The inverse of an $f \in \mathrm{Cont}(\partial)$ for the product $*$ is $-f$ (indeed, the inverse of $f\partial$ in the group of derivations is $-f\partial$), so we denote it so. The axioms **D1** and **D4** follow directly from the fact that $\preccurlyeq$ is a dominance relation on the Abelian group $(\mathrm{Cont}(\partial), +, 0)$.

Let $f, g \in \mathrm{Cont}(\partial)$. We have $f + g \preccurlyeq f$ or $f + g \preccurlyeq g$ by **D2** in $(\mathrm{Cont}(\partial), +, 0, \preccurlyeq)$. In view of (6.2) and Lemma 6.12, we deduce that **D2** holds in $(\mathrm{Cont}(\partial), *, 0, \preccurlyeq)$. We have $f * g * (-f) \sim g$ by [7, Lemma 1.3]. This implies that **D3** holds in $(\mathrm{Cont}(\partial), *, 0, \preccurlyeq)$ and, by Proposition 6.7, that this valued group is nearly Abelian. That **D5** holds follows from Lemma 6.14. □

Note that each non-trivial element $f$ in $\mathrm{Cont}(\partial)$ is scaling, and that the centraliser $\mathcal{C}(f) = C f$ is isomorphic to the underlying additive group of $C$. We also see since that group is divisible that $\mathrm{Cont}(\partial)$ is divisible.

**Proposition 6.16.** *The operation* $C \times \mathrm{Cont}(\partial) \longrightarrow \mathrm{Cont}(\partial)\,;\, (c, f) \mapsto c f$ *endows* $\mathrm{Cont}(\partial)$ *with a structure of* $C$*-group, and* $(\mathrm{Cont}(\partial), *, 0, (f \mapsto c f)_{c \in C}, \preccurlyeq)$ *is uniformly* $C$*-valued.*

**Proof.** The axioms are easily verified (e.g. that $g * (c f) * (-g) = c (g * f * (-g))$ for all $f$, $g \in \mathrm{Cont}(\partial)$ and $c \in C$ follows from the linearity of the Lie bracket). □

**Proposition 6.17.** *The valued group* $(\mathrm{Cont}(\partial), *, 0, \preccurlyeq)$ *is spherically complete.*

**Proof.** Let $\lambda$ be a non-zero limit ordinal and let $(f_\gamma)_{\gamma < \lambda}$ be pseudo-Cauchy in $\mathrm{Cont}(\partial)$. By [7, Corollary 1.1], we have

$$f * (-g) \sim f - g \tag{6.3}$$

for all $f, g \in \mathrm{Cont}(\partial)$. So $(f_\gamma)_{\gamma < \lambda}$ is pseudo-Cauchy in the valued field $C((\Gamma))$ in the sense of [4, p 64] (or pseudo-convergent as per Kaplansky [22]). We deduce by [23, p 193] that $(f_\gamma)_{\gamma < \lambda}$ pseudo-converges to an $f \in C((\Gamma))$. We have $f \preccurlyeq f_\gamma$ for a $\gamma < \lambda$ by [4, Lemma 2.2.3], so $f \in \mathrm{Cont}(\partial)$, and $f$ is a pseudo-limit of $(f_\gamma)_{\gamma < \lambda}$ in $\mathrm{Cont}(\partial)$ by definition and (6.3). We conclude with Proposition 5.27. □

Suppose now that $C$ is an ordered field. So $C((\Gamma))$ is an ordered field for the lexicographic ordering $f > 0 \Longleftrightarrow f \neq 0 \wedge f(\min \operatorname{supp} f) > 0$ (see [24, Section 1.5]) and the valuation is non-decreasing on $C((\Gamma))^{>}$. Therefore:

**Proposition 6.18.** *The structure* $(\mathrm{Cont}(\partial), *, 0, <, \preccurlyeq)$ *is an ordered c-valued group.*

Since $\mathrm{Ker}(\partial) = C$, the structure $(C((\Gamma)), +, \cdot, 0, 1, <, \preccurlyeq, \partial)$ is an H-field in the sense of [1, 2] if and only if $\partial(f) > 0$ for all $f > C$.

**Proposition 6.19.** *Suppose that* $(C((\Gamma)), \preccurlyeq, \partial)$ *is an H-field. Then the inverse group* $\big(\widecheck{\mathrm{Cont}}(\partial), \check{*}, 0, <\big)$ *is a growth order group whose finest c-dominance relation is* $\preccurlyeq$.

**Proof.** This is an ordered c-valued group by Remarks 3.4 and 4.5. In view of Lemma 6.14 and the definition of $<$, it satisfies **GOG1**, and $\preccurlyeq$ is the corresponding dominance relation. Each non-trivial element in $\mathrm{Cont}(\partial)$ is scaling in $(\mathrm{Cont}(\partial), *, 0, \preccurlyeq)$, so **GOG3** holds. Let us show that **GOG2** holds. Let $f, g \in \mathrm{Cont}(\partial)$ with $f, g > 0$ and $f > \mathcal{C}(g)$, i.e. $f \succ g$. We have $f \check{*} g - g \check{*} f = g * f - f * g \sim [g, f]$ by Lemma 6.12, so it suffices to show that $[g, f] > 0$, i.e. that $\frac{\partial(f)}{f} > \frac{\partial(g)}{g}$. Since $(C((\Gamma)), \preccurlyeq, \partial)$ is an H-field and $f \succ g$, this holds by [1, Lemma 1.4]. □

**Remark 6.20.** The relation $f \succ g \Longrightarrow \frac{\partial(f)}{f} > \frac{\partial(g)}{g}$ is equivalent to being an H-field modulo our assumptions on $\partial$. So we have an equivalence between two seemingly unrelated classes of objects: growth order groups and H-fields.

### 6.3 Nearly Abelian quotients of ordered valued groups

Let $(\mathcal{G}, \cdot, 1, <, \preccurlyeq)$ be an ordered valued group satisfying **GA**. Let us see how $\mathcal{G}$ can be decomposed into nearly Abelian parts. Recall by Proposition 4.28 that $(\mathcal{G}, \cdot, 1, <, \underline{\prec\!\!\prec})$ is an ordered valued group satisfying **GA** and that $\underline{\prec\!\!\prec}$ is coarser than $\preccurlyeq$. We write val* for the valuation on $\mathcal{G}$ corresponding to $\underline{\prec\!\!\prec}$. Since $\preccurlyeq$ is convex on $\mathcal{G}^{\geqslant}$, so is $\underline{\prec\!\!\prec}$. Therefore, the quotient set $\mathcal{G}/\asymp\!\!\asymp$ is totally ordered by the relation

$$\mathrm{val}^*(f) < \mathrm{val}^*(g) \Longleftrightarrow f \prec\!\!\prec g \Longleftrightarrow \mathrm{val}^*(f) \cap \mathcal{G}^{\geqslant} < \mathrm{val}^{\times}(g) \cap \mathcal{G}^{\geqslant}.$$

Let $\beta = \mathrm{val}^*(f_0) \in \mathrm{val}^*(\mathcal{G})$ and set

$$\mathcal{G}_{\leqslant\beta} := \{g \in \mathcal{G} : g \,\underline{\prec\!\!\prec}\, f_0\} \qquad \text{and} \qquad \mathcal{G}_{<\beta} := \{g \in \mathcal{G} : g \prec\!\!\prec f_0\}.$$

By Proposition 4.28, those are $\preccurlyeq$-initial subsets of $\mathcal{G}$. We deduce (see Section 4.1) that they are convex subgroups of $\mathcal{G}$. The validity of **D2**, **D3** and **D4** for $\underline{\prec\!\!\prec}$ entails that $\mathcal{G}_{<\beta}$ is normal subgroup of $\mathcal{G}_{\leqslant\beta}$. We deduce:

**Proposition 6.21.** *Suppose that for all $f, g \asymp\!\!\asymp f_0$, we have $[f, g] \prec\!\!\prec f \Longrightarrow f \asymp g$. Then the quotient $\mathcal{G}_{\leqslant\beta}/\mathcal{G}_{<\beta}$, together with the ordering and valuation of Proposition 4.13, is a nearly Abelian ordered c-valued group satisfying* **GA**.

**Remark 6.22.** This result is significant, for it points to a representation of growth order groups as subgroups of infinite semidirect products of nearly Abelian quotients $\mathcal{G}_{\leqslant\beta}/\mathcal{G}_{<\beta}$ where $\beta$ ranges in val*$(\mathcal{G})$. This applies for certain groups $\mathcal{P}$ of parabolic elements in an H-field $K$ with composition and inversion (see Section 4.1). In future work, we will show that this applies to the field of finitely nested hyperseries of Example 4.16, give synthetic conditions on $\mathcal{G}$ under which Proposition 6.21 applies, give conditions for the existence of such sections, and construct such infinite semidirect products.

## 7 Regular equations over nearly Abelian valued groups

In this section, we fix a ring $A$ of characteristic 0, and a nearly Abelian $A$-valued $A$-group $(\mathcal{G}, \cdot, 1, (g \mapsto g^a)_{a \in A}, \preccurlyeq)$ which is $A$-torsion-free. We recall (see Remark 3.11) that each Abelian group $\mathcal{C}_\rho$ for $\rho \in \mathrm{val}(\mathcal{G})$ is an $A$-module for the operation $a\,\mathrm{res}(f) := \mathrm{res}(f^a)$ for $a \in A$ and $\mathrm{res}(f) \in \mathcal{C}_\rho$. Since $\mathcal{G}$ is $A$-torsion-free, this module has no $A$-torsion. Since $A$ has characteristic 0, the group $\mathcal{G}$ is in particular $\mathbb{Z}$-valued and torsion-free.

### 7.1 Representation of unary terms

Consider a multiplicative copy $(y^A, \cdot, 1)$ of the underlying additive group of $A$. So $y^A$ is Abelian and torsion-free. Note that $y^A$ is an $A$-module for the operation $(a, y^b) \mapsto y^{ab}$. We have [26, Theorem 5] a free product of $A$-groups $\mathcal{G} *^A y^A$. In view of Proposition 6.9, both $\mathcal{G}$ and $y^A$ are CSA*-groups as per [27, Definition 13], so $\mathcal{G} *^A y^A$ can be obtained as the reunion of subgroups $\mathcal{G} *^A y^A = \bigcup_{m \in \mathbb{N}} \mathcal{T}_m$ where $\mathcal{T}_0 = \mathcal{G} *^{\mathbb{Z}} y^{\mathbb{Z}}$ and

$$\mathcal{T}_{m+1} := \{t_1^{a_1} \cdots t_n^{a_n} : n \in \mathbb{N} \wedge t_1, \ldots, t_n \in \mathcal{T}_m \wedge a_1, \ldots, a_n \in A\}$$

for all $m \in \mathbb{N}$. If $\mathcal{H} \supseteq \mathcal{G}$ is an $A$-group extension and $h \in \mathcal{H}$, then by the universal property of the free product, there is a unique homomorphism of $A$-groups $\mathrm{ev}_h : \mathcal{G} *^A y^A \longrightarrow \mathcal{H}$ that fixes $\mathcal{G}$ pointwise and sends $y$ to $h$. We write $t(h) := \mathrm{ev}_h(t)$ for each $t \in \mathcal{G} *^A y^A$. Note that we may apply this to $\mathcal{H} = \mathcal{G} *^A y^A$ itself.

We denote by $\mathcal{L}_{\mathrm{Ag}}(\mathcal{G})$ the first-order language of $A$-groups expanded with constants for elements in $\mathcal{G}$. By induction on the complexity of terms, we see that for each unary term $t$ in $\mathcal{L}_{\mathrm{Ag}}(\mathcal{G})$, there is a canonical element $\hat{t} \in \mathcal{G} *^A y^A$ such that the evaluation of $t$ at an $h \in \mathcal{H}$ in any $A$-group extension is $\hat{t}(h)$.

Again by the universal property of free products, there is a unique $A$-group homomorphism $\alpha : \mathcal{G} *^A y^A \longrightarrow A$, such that $\alpha(g) = 0$ for all $g \in \mathcal{G}$ and $\alpha(y) = 1$. We say that a $t \in \mathcal{G} *^A y^A$ is *regular* if $\alpha(t) \neq 0$.

Given $t \in \mathcal{G} *^A y^A$ and $f \in \mathcal{G}$, we write $t_f := t(f)^{-1} t(fy)$, so $t(fy) = t(f)\, t_f$ and $t_f(1) = 1$. The map $t \mapsto \alpha(t(fy))$ satisfies the defining properties of $\alpha$, so it coincides with $\alpha$. We thus have $\alpha(t_f) = \alpha(t(fy)) = \alpha(t)$.

## 7.2 Solving regular equations

**Proposition 7.1.** *Let $t \in \mathcal{G} *^A y^A$ with $t(1) = 1$. For all $f \in \mathcal{G}^{\neq}$, we have $\mathrm{res}(t(f)) = \alpha(t)\,\mathrm{res}(f)$ if $t$ is regular, and $t(f) \prec f$ otherwise.*

**Proof.** We prove this by induction on $m \in \mathbb{N}$ where $t \in \mathcal{T}_m$. If $m = 0$, then we can write $t = g_1\, y^{\alpha_1} \cdots g_l\, y^{\alpha_l}$ where $l > 0$, $g_1, \ldots, g_l \in \mathcal{G}$ and $\alpha_1, \ldots, \alpha_l \in \mathbb{Z}$. We have $\alpha(t) = \alpha_1 + \cdots + \alpha_l$ and $t(1) = g_1 \cdots g_l$. Since $t(1) = 1$, we can rewrite $t = \prod_{i=1}^{n} \tilde{t_i}$ where

$$\tilde{t_i} := \left( \prod_{j=i+1}^{l} g_j \right)^{-1} y^{\alpha_i} \left( \prod_{j=i+1}^{l} g_j \right)$$

for each $i \in \{1, \ldots, n\}$. Now $\mathrm{res}(\tilde{t_i}(f)) = \mathrm{res}(f^{\alpha_i}) = \alpha_i\, \mathrm{res}(f)$ for all $i \in \{1, \ldots, n\}$ by Proposition 6.7. Recall that $\mathcal{C}_{\mathrm{val}(f)}$ is $A$-torsion-free, hence torsion-free. We conclude with Proposition 2.31.

Suppose now that $m > 0$ and that the result holds on $\mathcal{T}_{m-1}$,. Write $t = t_1^{a_1} \cdots t_n^{a_n}$ accordingly, where $a_1, \ldots, a_n \in A$ and $t_1, \ldots, t_n \in \mathcal{T}_{m-1}$. Note that $\alpha(t) = \sum_{i=1}^{n} a_i\, \alpha(t_i)$. Since $t(1) = t_1(1)^{a_1} \cdots t_n(1)^{a_n} = 1$, we may rewrite $t = \prod_{i=1}^{n} \hat{t_i}$ where

$$\hat{t_i} := \left( \prod_{j=1}^{i-1} t_j(1)^{a_j} \right) t_i^{a_i} \left( \prod_{j=1}^{i} t_j(1)^{a_j} \right)^{-1}$$

for each $i \in \{1,\ldots,n\}$. We have $\hat{t_i}(1) = 1$ for each $i \in \{1,\ldots,n\}$ by construction. Note since $\mathcal{T}_m$ is a subgroup containing $\mathcal{G}$ that each $\hat{t_i}$ is in $\mathcal{T}_{m-1}$. Set $J := \{i \in \{1,\ldots,n\} : a_i\, \alpha(t_i) \neq 0\}$. For $j \in J$, we have $\mathrm{res}(\hat{t_i}(f)) = a_i\, \alpha(t_i)\, \mathrm{res}(f)$ by the induction hypothesis. For $j \in \{1,\ldots,n\} \setminus J$, we have $\hat{t_i}(f) \prec f$ by the induction hypothesis. If $\alpha(t) \neq 0$, then $J$ must be non-empty since $A$ is torsion-free. Furthermore $\mathcal{C}_{\mathrm{val}(f)}$ has no $A$-torsion, so Proposition 2.31 gives $\mathrm{res}(t(f)) = \sum_{j \in J} \mathrm{res}(\hat{t_j}(f)) = \sum_{j \in J} a_i\, \alpha(t_i)\, \mathrm{res}(f) = \alpha(t)\, \mathrm{res}(f)$. If $\alpha(t) = 0$, then we either have $J = \varnothing$ and thus $t(f) \prec f$ by Proposition 2.4, or $J \neq 0$ and $t(f) \prec f$ by Proposition 2.31. This concludes the inductive proof. □

**Corollary 7.2.** *Let $\lambda$ be a non-zero limit ordinal, let $(g_\gamma)_{\gamma<\lambda}$ be a pseudo-Cauchy sequence in $\mathcal{G}$ and let $t \in \mathcal{G} *^A y^A$ be regular. Then $(t(g_\gamma))_{\gamma<\lambda}$ is pseudo-Cauchy. If moreover $(g_\gamma)_{\gamma<\lambda}$ pseudo-converges to a $g \in \mathcal{G}$, then $(t(g_\gamma))_{\gamma<\lambda}$ pseudo-converges to $t(g)$.*

**Proof.** By Proposition 5.26, we may assume that $(g_\gamma)_{\gamma<\lambda}$ has a pseudo-limit $g\in\mathcal{G}$. By Lemma 5.24, we may assume that $g=1$ and that $t(g)=1$. We have $\operatorname{res}(t(g_\gamma))=\alpha(t)\operatorname{res}(g_\gamma)$ for all $\gamma<\lambda$ by Proposition 7.1. We deduce that $(t(g_\gamma))_{\gamma<\lambda}$ pseudo-converges to 1. □

**Corollary 7.3.** *Suppose that $A=\mathbb{Z}$ and that $\mathcal{G}$ is ordered c-valued. Let $t\in\mathcal{G}*y^{\mathbb{Z}}$ be regular with $t(1)=1$. Then for all $f\in\mathcal{G}$, we have $t(f)\geqslant 1\Longleftrightarrow f^{\alpha(t)}\geqslant 1$.*

**Lemma 7.4.** *Let $t\in\mathcal{G}*^A y^A$ be regular with $t(1)\neq 1$. Then for all $f\in\mathcal{G}$, we have*

$$\operatorname{val}(t(f))<\operatorname{val}(t(1))\Longleftrightarrow\alpha(t)\operatorname{res}(f)=-\operatorname{res}(t(1)).$$

**Proof.** We have $t_1(1)=1$ and $\alpha(t_1)=\alpha(t)\neq 0$, so $\operatorname{res}(t_1(f))=\alpha(t)\operatorname{res}(f)$ by Proposition 7.1, since $\alpha(t)\operatorname{res}(f)\neq 0$ since $\mathcal{C}_\rho$ is $A$-torsion-free. Since $t(1)\neq 1$, it follows by Proposition 2.31 that $\operatorname{val}(t(f))<\operatorname{val}(t(1))$ if and only if $\alpha(t)\operatorname{res}(f)=-\operatorname{res}(t(1))$. □

Given $t\in\mathcal{G}*^A y^A$, we write $\operatorname{Val}(t):=\{\rho\in\operatorname{val}(\mathcal{G}):\exists g\in\mathcal{G},\operatorname{val}(t(g))\leqslant\rho\}$.

**Lemma 7.5.** *Let $t\in\mathcal{G}*^A y^A$ be regular and set $\alpha:=\alpha(t)$. There is a pseudo-Cauchy sequence $(f_\gamma)_{\gamma<\lambda}$ for some $\lambda\in\mathbf{On}$ such that one of the following occurs:*

*a) $\lambda=\gamma+1$ is a successor and*

$$\begin{aligned}\forall\eta<\lambda,\operatorname{res}(t(f_\eta))\in\alpha\,\mathcal{C}_{\operatorname{val}(t(f_\eta))}\quad&\Longleftrightarrow\quad\eta\neq\gamma,\\ \operatorname{Val}(t)\quad&=\quad[\operatorname{val}(t(f_\gamma)),+\infty).\end{aligned}$$

*b) $\lambda=\gamma+1$ is a successor and*

$$\begin{aligned}\forall\eta<\lambda,\operatorname{res}(t(f_\eta))&\in\alpha\,\mathcal{C}_{\operatorname{val}(t(f_\eta))},\\ \operatorname{Val}(t)&=\operatorname{val}(\mathcal{G}),\\ t(f_\gamma)&=1,\\ \forall g\in\mathcal{G}\setminus\{f_\gamma\},t(g)&\neq 1.\end{aligned}$$

*c) $\lambda$ is a non-zero limit and $(f_\gamma)_{\gamma<\lambda}$ has no pseudo-limit in $\mathcal{G}$ and*

$$\begin{aligned}\forall\eta<\lambda,\operatorname{res}(t(f_\eta))&\in\alpha\,\mathcal{C}_{\operatorname{val}(t(f_\eta))},\\ \operatorname{Val}(t)&=\bigcup_{\eta<\lambda}[\operatorname{val}(t(f_\eta)),+\infty).\end{aligned}$$

**Proof.** We will construct $(f_\gamma)_{\gamma\leqslant\lambda}$ by induction, insuring that the sequence $(\operatorname{val}(t(f_\gamma y)))_{\gamma<\lambda}$ is strictly decreasing. We start with $f_0:=1$. Let $\eta$ be such that $(f_\gamma)_{\gamma<\eta}$ is already defined and that each initial segment of that sequence satisfies the conditions.

Suppose that $\eta$ is a limit. If $(f_\gamma)_{\gamma<\eta}$ has a pseudo-limit $f$ in $\mathcal{G}$, then by Corollary 7.2 we may extend the sequence with $f_\eta:=f$. If $(f_\gamma)_{\gamma<\eta}$ has no pseudo-limit in $\mathcal{G}$, by Lemma 7.4 there is no $f\in\mathcal{G}$ such that $\operatorname{val}(t(f_\gamma\,(f_\gamma^{-1}\,f)))<\operatorname{val}(t(f_\gamma))$ for all $\gamma<\eta$. Therefore, setting $\lambda:=\eta$, we see that c holds.

Suppose now that $\eta=\gamma+1$ is a successor. If $t(f_\gamma)=1$, then $t(f_\gamma y)$ satisfies the conditions of Corollary 7.3, so $f_\gamma$ is the unique solution of $t(y)=1$, thus we may set $\lambda=\eta$ and find that b holds. We next assume that $t(f_\gamma)\neq 1$. Suppose furthermore that $\operatorname{res}(t(f_\gamma))\in\alpha\,\mathcal{C}_{\operatorname{val}(t(f_\gamma))}$. Then pick a $g\in\mathcal{G}$ with $\alpha\operatorname{res}(f)=-\operatorname{res}(t(f_\gamma))$ and set $f_\eta:=f_\gamma\,g$. We have $\operatorname{val}(t(f_\eta))=\operatorname{val}(t(f_\gamma y)(g))<\operatorname{val}(t(f_\gamma))$ by Lemma 7.4. Assume for contradiction that there is a $\mu\leqslant\gamma$ with $g\succcurlyeq f_\mu$. We have

$$\operatorname{val}(t(f_{\mu+1}))\succ\operatorname{val}(t(f_\mu))=\operatorname{val}(t(f_\mu\,(f_\mu^{-1}\,f_\gamma g))),$$

so Lemma 7.4 implies that $\alpha\,\mathrm{res}(f_\mu\,(f_\mu^{-1}\,f_\gamma\,g))=-\mathrm{res}(t(f_\mu))=-\alpha\,\mathrm{res}(f_\mu)$. For $\mu<\eta$, we have $t(f_\mu\,y)(f_\mu^{-1}\,f_\gamma\,g)=t(f_\eta)\prec t(f_\mu\,y)$. We deduce by Lemma 7.4 that we must have $\mathrm{res}(f_\mu^{-1}\,f_\gamma\,g)=\mathrm{res}(f_\mu^{-1}\,f_\gamma)$. In view of Proposition 2.31, this entails that $g\prec f_\mu\,f_\gamma^{-1}$. We deduce with Lemma 5.25 that $(f_\mu)_{\mu<\eta+1}$ is pseudo-Cauchy. Lemma 7.4 entails that $(f_\mu)_{\mu<\eta+1}$ is pseudo-Cauchy satisfies the condition.

Lastly, if $\mathrm{res}(t(f_\gamma))\notin\alpha\,\mathcal{C}_{\mathrm{val}(t(f_\gamma))}$, then by Lemma 7.4, we have $V(t)=[\mathrm{val}(t(f_\gamma)),+\infty)$. So we may set $\lambda:=\eta$ and see that a holds. This concludes the inductive construction. Since $(\mathrm{val}(t(f_\gamma)))_{\gamma<\lambda}$ is strictly decreasing for each $\lambda$, there must be an ordinal $\lambda$ such that this process stops, i.e. for which one of the cases occurs. □

**Corollary 7.6.** *If $t\in\mathcal{G}*^A y^A$ is regular, then there is at most one $f\in\mathcal{G}$ with $t(f)=1$.*

**Theorem 7.7.** *Assume that $\mathcal{G}$ is spherically complete. Let $t\in\mathcal{G}*^A y^A$ be regular and suppose that $\alpha(t)\,\mathcal{C}_\rho=\mathcal{C}_\rho$ for all $\rho\in\mathrm{val}(\mathcal{G})$. Then there is a unique $\psi\in\mathcal{G}$ with $t(\psi)=1$.*

**Proof.** The cases a and c in Lemma 7.5 cannot occur, so Lemma 7.5(b) occurs. □

**Corollary 7.8.** *Assume that $\mathcal{G}$ is spherically complete and that each group $\mathcal{C}_\rho$, $\rho\in\mathrm{val}(\mathcal{G})$ is divisible. Then $\mathcal{G}$ is divisible.*

**Theorem 7.9.** *Suppose that $A=\mathbb{Z}$ and that $\mathcal{G}$ is ordered valued. Let $t\in\mathcal{G}*y^{\mathbb{Z}}$ be regular. Then the function $f\mapsto t(f)$ is strictly increasing if $\alpha(t)>0$, and strictly decreasing if $\alpha(t)<0$.*

**Proof.** Considering $t^{\alpha(t)}$ if necessary, we may assume that $\alpha(t)>0$. Since we can multiply $t$ by elements in $\mathcal{G}$ while preserving the sign of $\alpha(t)$, it suffices to show that the set of elements $f\in\mathcal{G}$ with $t(f)<1$ is an initial subset of $(\mathcal{G},<)$.

Let $(f_\gamma)_{\gamma\leqslant\lambda}$ be as in Lemma 7.5 with respect to $t$. If Lemma 7.5(a) occurs and $\lambda=\gamma+1$, then the fact that $\mathrm{res}(t(f_\gamma))\notin\alpha(t)\,\mathcal{C}_{\mathrm{val}(t(f_\gamma))}$ implies that for all $f\in\mathcal{G}$, we have $f_\gamma^{-1}\,f\lessdot t(f_\gamma)$ or $f_\gamma^{-1}\,f\gtrdot t(f_\gamma)$. In that first case, we have $t(f)<t(1)$, and in the second case, we have $t(f)>1$, hence the result. If Lemma 7.5(b) occurs and $\lambda=\gamma+1$, then Corollary 7.3 for $t(f_\gamma\,y)$ yields the result. Lastly if Lemma 7.5(c) occurs, then the adjunction of $f$ does not extend $(f_\gamma)_{\gamma<\lambda}$ into a pseudo-Cauchy sequence. This means by Lemma 5.25 that $(f_\gamma^{-1}\,f)\sim f_\gamma^{-1}\,f_{\gamma+1}$ for some $\gamma<\lambda$. But then Lemma 7.4 for $t(f_\gamma y)$ entails that that $t(f)<1$ if $f_\gamma^{-1}\,f\lessdot f_\gamma^{-1}\,f_{\gamma+1}$ and that $t(f)>1$ if $f_\gamma^{-1}\,f\gtrdot f_\gamma^{-1}\,f_{\gamma+1}$, hence the result. □

**Theorem 7.10.** *Let $\mathcal{H}$ be a residually $A$-nilpotent $A$-group whose quotient modules are $A$-torsion-free. Then any equation $t(y)=1$ for $t\in\mathcal{H}*^A y^A$ with $\alpha(t)\in A^\times$ has a unique solution in a residually $A$-nilpotent $A$-group whose quotient modules are $A$-torsion-free. If $\mathcal{H}$ is $A$-nilpotent, then the solution lies in $\mathcal{H}$ itself.*

**Proof.** We endow $\mathcal{H}$ with the lower central valuation, making it an $A$-valued $A$-group with value set $(\beta,>)$ for a $\beta\leqslant\omega$ (see Proposition 3.6). Proposition 5.14 entails that the completion $\tilde{\mathcal{H}}$ satisfies the same properties. It is nearly Abelian by Proposition 6.6. Since $v(\mathcal{H})=(\beta,>)$, the valued-group $\tilde{\mathcal{H}}$ is spherically complete. We have $\mathcal{C}_\rho=\alpha(t)\,(\alpha(t)^{-1}\,\mathcal{C}_\rho)\subseteq\alpha(t)\,\mathcal{C}_\rho$ for all $\rho\in v(\tilde{\mathcal{H}})$, so we conclude by applying Theorem 7.7 to $\tilde{\mathcal{H}}$. If $\mathcal{H}$ is $A$-nilpotent, then it is complete, so $\tilde{\mathcal{H}}=\mathcal{H}$. □

**Acknowledgments.** We thank Mickaël Matusinski for our helpful conversations with him on this topic. We thank Simon André and Andrei Jaikin-Zapirain for their answers to our questions, and Sylvy Anscombe for her precious advice.